\documentclass[11pt]{article}
\usepackage{amsmath,amsfonts,amssymb,amscd,amsthm,amsbsy,epsf}
\usepackage{graphicx}

\newtheorem{prop}{Proposition}[section]
\newtheorem{theorem}[prop]{Theorem}
\newtheorem{lemma}[prop]{Lemma}

\newtheorem{corollary}[prop]{Corollary}
\newtheorem{remark}[prop]{Remark}

\pdfoutput=1

\begin{document}

\title{Lifted Heegaard Surfaces and Virtually Haken Manifolds}

\author{Yu Zhang\thanks{E-mail: yz26@buffalo.edu}\\ University at Buffalo, The State University of New York}

\date{ }

\maketitle

\begin{abstract}
  In this paper, we give  infinitely  many  non-Haken  hyperbolic  genus three $3$-manifolds
  each of which has a finite cover whose induced Heegaard surface   from  some
  genus three Heegaard surface  of the base manifold  is   reducible
   but can be compressed into an incompressible surface.
  This result supplements \cite{cg} and extends \cite{mmz}.
\end{abstract}

\section{Introduction}

It was shown in \cite{cg} that if a Heegaard splitting of an irreducible closed $3$-manifold $M$ is
weakly reducible then either the Heegaard splitting is reducible or $M$ contains an incompressible
surface of positive genus. This  result  motivates an approach to the well known virtual Haken
conjecture which, with the current knowledge,  is reduced to the following conjecture:
 every closed hyperbolic $3$-manifold is virtually Haken, i.e.
 has a finite cover which is a Haken $3$-manifold.
 That is, to prove that a given closed hyperbolic $3$-manifold  is virtually Haken, it suffices to find a
  finite cover which has an irreducible but weakly reducible
 Heegaard splitting.
In \cite{mmz},  families of non-Haken but virtually Haken hyperbolic $3$-manifolds were found using this
approach. These manifolds were obtained by Dehn surgeries on some $2$-bridge knots in $S^3$ and thus are
genus two $3$-manifolds. In fact it was showed there that each of these manifolds has a finite cover
whose induced Heegaard surface from some genus two Heegaard surface  of the base manifold is weakly
reducible and can be compressed into an incompressible surface, without the need to know whether the
Heegaard surface of the cover is irreducible or not (we suspect that it is irreducible).

The main purpose of this paper is to illustrate two points concerning the above works. One point is to
show that the method used in \cite{mmz} can be generalized to find  an infinite  family of closed
non-Haken but virtually Haken hyperbolic genus three $3$-manifolds. The other point is to show that each
manifold of our family has a finite cover whose induced  Heegaard surface from some genus three Heegaard
surface of the base manifold is actually reducible but can still be compressed into an incompressible
surface, which is a phenomenon supplementing \cite{cg}. Our manifolds are obtained by Dehn surgeries on
some pretzel knots in $S^3$.

Let $K=(p,\pm 3, q)$ be a pretzel knot in $S^3$ with $p, q$ odd and $|p|, |q|\geqslant 3$, and let
$M_K=S^3\setminus\stackrel{\circ}{N}(K)$ be the exterior of $K$. Let $M_K^3$ be the $3$-fold cyclic cover
of $M_K$. We give $\partial M_K$ the standard meridian-longitude coordinates and $\partial M_K^3$ the
induced meridian-longitude coordinates. So a slope in such a torus can be identified with  a rational
number $m/n$ where  $m$ is the meridian coordinate  and $n$ the longitude coordinate. By \cite{la}, $K$
is a tunnel number two knot and thus $M_K$ is a genus three manifold.

\begin{theorem}\label{theorem 1}
For the pretzel knot $K=(p,\pm 3, q)$, the induced Heegaard surface of $M_K^3$ from some genus three
Heegaard surface of $M_K$ is reducible and can be compressed into an essential surface $S$ in $M_K^3$.
Moreover, $S$ remains essential in  every Dehn filling of  $M_K^3$ with slope $m/n, (m, n)=1,
|m|\geqslant 2$. Thus every Dehn filling of  $M_K$ with slope $3m/n, (3m, n)=1, |m|\geqslant 2$, yields
a virtually Haken $3$-manifold.
\end{theorem}

We now explain how the results described in the second paragraph of this section follow from Theorem
\ref{theorem 1}. As we have noted, $M_K$ is a genus three manifold. Thus every Dehn filling of $M_K$ is
of genus at most three. As $K$ is a hyperbolic small knot by \cite{o} (here small means no closed
embedded essential surfaces in $M_K$), it follows from \cite{rs} that except for finitely many lines in
the Dehn filling plane of $M_K$, all remaining Dehn fillings of $M_K$ are genus three manifolds, which
we may also assume to be \newline (1) hyperbolic, by Thurston's hyperbolic Dehn surgery theorem, and
\newline (2) non-Haken, by \cite{h}.\newline
 Hence infinitely many of $M_K(3m/n)$ given in
Theorem \ref{theorem 1} are genus three non-Haken hyperbolic $3$-manifolds. Finally we just need to note
that $M_K(3m/n)$ is covered by $M_K^3(m/n)$ and that each Heegaard splitting of $M_K$ induces a Heegaard
splitting on $M_K(3m/n)$.

 The proof of Theorem \ref{theorem 1} is given in Section 3, after
 some preliminary preparations in Section 2.
Using a  similar method, we shall also give  a new  proof of \cite[Corollary 4(b)]{o}
 in case of pretzel knots. This is the content of Section 4.

\section{Preliminary}\label{pre}

\textbf{Heegaard Splittings.} A \textit{Heegaard splitting} $M=W_1\cup_F W_2$ of a compact $3$-manifold $M$ is a decomposition of $M$ into two compression bodies $W_1$ and $W_2$ with common positive boundary $F$. A Heegaard splitting $M=W_1\cup_F W_2$ is \textit{reducible} if there exist essential disks $(D_1, \partial D_1)\subset (W_1, F)$ and $(D_2, \partial D_2)\subset (W_2, F)$ such that $\partial D_1=\partial D_2$. Otherwise, it is \textit{irreducible}. If neither $W_1$ nor $W_2$ is trivial and there do not exist essential disks $(D_1, \partial D_1)\subset (W_1, F)$ and $(D_2, \partial D_2)\subset (W_2, F)$ such that $\partial D_1\cap \partial D_2=\emptyset$, then the Heegaard splitting is \textit{strongly irreducible}. Otherwise, it is \textit{weakly reducible}. We call a Heegaard splitting $M=W_1\cup_F W_2$ \textit{stabilized} if there exist essential disks $(D_1, \partial D_1)\subset (W_1, F)$ and $(D_2, \partial D_2)\subset (W_2, F)$ such that $\partial D_1$ and $\partial D_2$ intersect at a single point. It is known that every reducible splitting of an irreducible manifold is stabilized.

\begin{figure}
\begin{center}
\includegraphics[width=4in]{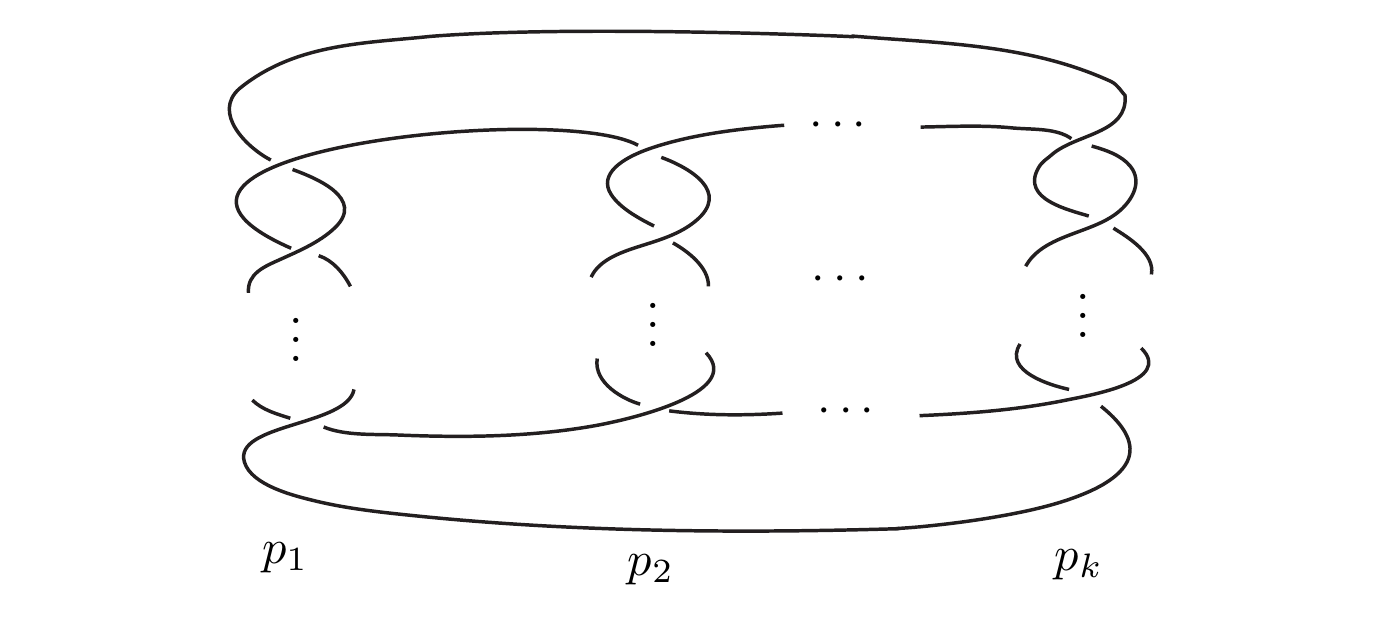}
\end{center}
\caption{\label{pretzel} Pretzel link}
\end{figure}

\noindent\textbf{Pretzel Links.} A \textit{pretzel link} is a special kind of link. A pretzel link which is also a knot is a \textit{pretzel knot}. In the standard projection of the $(p_1, p_2, \cdots, p_k)$-pretzel link, there are $p_i$ left-handed crossings in the $i$th tangle, see Figure \ref{pretzel}. Obviously, the $(p_1, p_2, \cdots, p_k)$-pretzel link is link-equivalent to the $(p_i, p_{i+1}, \cdots, p_k, p_1, \cdots, p_{i-1})$-pretzel link.

For a pretzel link $K=(p_1, p_2, \cdots, p_k)$, the number of components $|K|$ of $K$ is given by
\begin{equation}\label{components}
|K|=\begin{cases}
1\ \ \ \ \ \ \ \ \ \ \ \ \ \ \ \ \ \ \ \ \ \ \ \ \ \ \ \text{if each $p_i$ is odd and $k$ is odd}\\
2\ \ \ \ \ \ \ \ \ \ \ \ \ \ \ \ \ \ \ \ \ \ \ \ \ \ \ \text{if each $p_i$ is odd and $k$ is even}\\
\#\{i: \text{$p_i$ is even}\}\ \ \ \ \ \ \ \text{if some $p_i$'s are even.}
\end{cases}
\end{equation}

\noindent \textbf{The Whitehead Graphs.} Let $H$ be a handlebody. A finite set of pairwise  disjoint
simple closed curves $A$ in $\partial H$ is said to be {\it  separable} in $H$  if and only if $\partial
H- A$ is compressible in $H$. Let $D$ be a compression disk system of $H$, i.e., a set of embedded
essential disks in $H$ which  compress $H$ into a $3$-ball $B$. Each disk $d$ in $D$ has two copies $d^+,
d^-$ in $\partial B$. The {\it Whitehead graph of $A\subset \partial H$ with respect to $D$ }, denoted by $WG(D,A)$, is the
graph on $\partial B$ taking   $\cup\; d^\pm$ as vertices and taking  the line segments $[\partial B-
int(\cup\; d^\pm)]\cap A$ as edges.
 It is shown in
\cite{s} that if the Whitehead graph  is connected and has no cut vertex, then $A$ is non-separable in
$H$,  and if the graph is disconnected, then $A$ must be separable. When the graph has a cut vertex $v$,
a {\it Whitehead  automorphism} corresponding to $v$ can be made  to transform the graph into an
equivalent graph (by changing the disk system $D$) which has less complexity (i.e. the number of edges).
So after a finitely many Whitehead automorphisms, we may end  up with a disconnected graph or a
connected graph with no cut vertices. We refer to \cite{s} for details about how to make the Whitehead
automorphism at a cut vertex. The following two elementary lemma and corollary will be handy  in the proof of Theorem
\ref{theorem 1}.

\begin{figure}
\begin{center}
\includegraphics[width=4in]{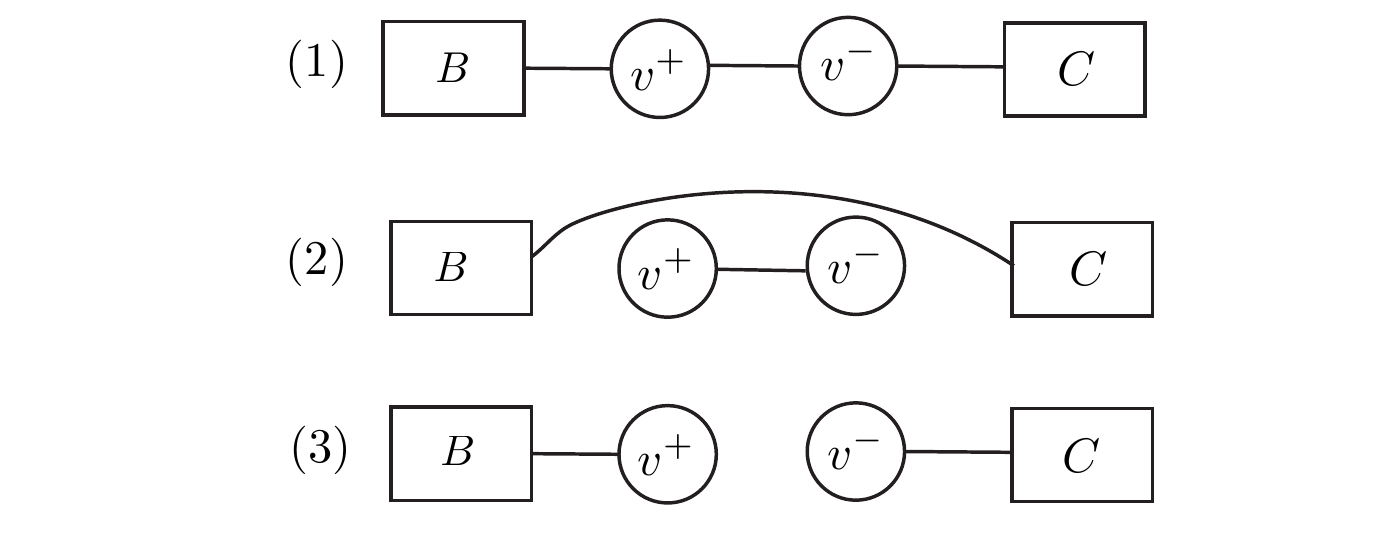}
\end{center}
\caption{\label{lem2}A graph contains a path with $v^+$ and $v^-$.}
\end{figure}

\begin{lemma}\label{lemma 2}
If the Whitehead graph of $A\subset \partial H$ consists of
 as two subgraphs $B$ and $C$ connected by a path
with two vertices $\{ v^+, v^-\}$ (as shown in part (1) of Figure \ref{lem2}), then $A$ is separable in
$H$.
\end{lemma}

\textbf{Proof:} In  Figure \ref{lem2} (1), the vertex  $v^+$ is a cut vertex. Applying the Whitehead
automorphism to $v^+$, we get a new graph which looks like Figure \ref{lem2} (2) or (3), both being
disconnected. So
$A$ is separable.

\begin{corollary}\label{lemma 1}
If the Whitehead graph of $A\subset \partial H$ has a vertex of valence one, then $A$ is separable in
$H$.
\end{corollary}

Lastly in this section  we record the Multi-Handle Addition Theorem  given in \cite{le}.

\begin{theorem}\label{mha}{\rm\cite{le}}
 Let $C=\{c_1, \cdots, c_n\}$ be a set of  pairwise disjoint simple closed curves in
the boundary of a  handlebody $H$ of genus $k>0$. If the following conditions are satisfied:
\newline
(0) $\partial H-C$ is incompressible in $H$,
\newline
(1) for each $j$, $\partial H-(C-c_j)$ is compressible in $H$, i.e., $C-c_j$ does not bind the free group $F_k$,
\newline
(p) for any $(n-p)$-element subfamily $C'$ of $C$, $C'$ does not bind any free factor $F_{k-p+1}$ of $F_k$,
\newline
 (n-1) for any $c_j\in C$, $c_j$ does not bind a free factor $F_{k-n+2}$ of $F_k$.
\newline
Then the 3-manifold obtained by adding $n$ 2-handles to $H$ along $C$ has incompressible boundary.
\end{theorem}

See \cite{le} for  the term ``bind a free factor''.

\section{Proof of Theorem \ref{theorem 1}}

\begin{figure}
\begin{center}
\includegraphics[width=4in]{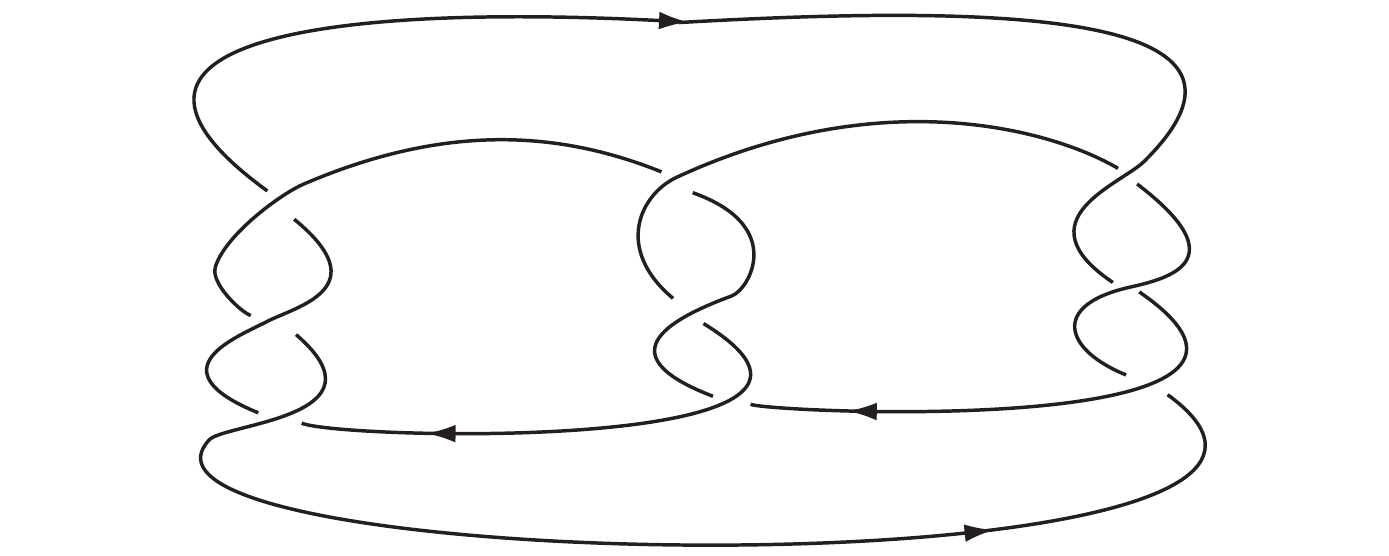}
\end{center}
\caption{\label{333} The standard diagram of $(3, 3, 3)$-pretzel knot}
\end{figure}

We first give  a detailed proof  when $K$ is the  $(3, 3, 3)$-pretzel knot (Figure \ref{333} shows its
standard diagram) and then indicate how to extend the proof  to work for general $K=(p,\pm3,q)$.

\begin{figure}
\begin{center}
\includegraphics[width=4.5in]{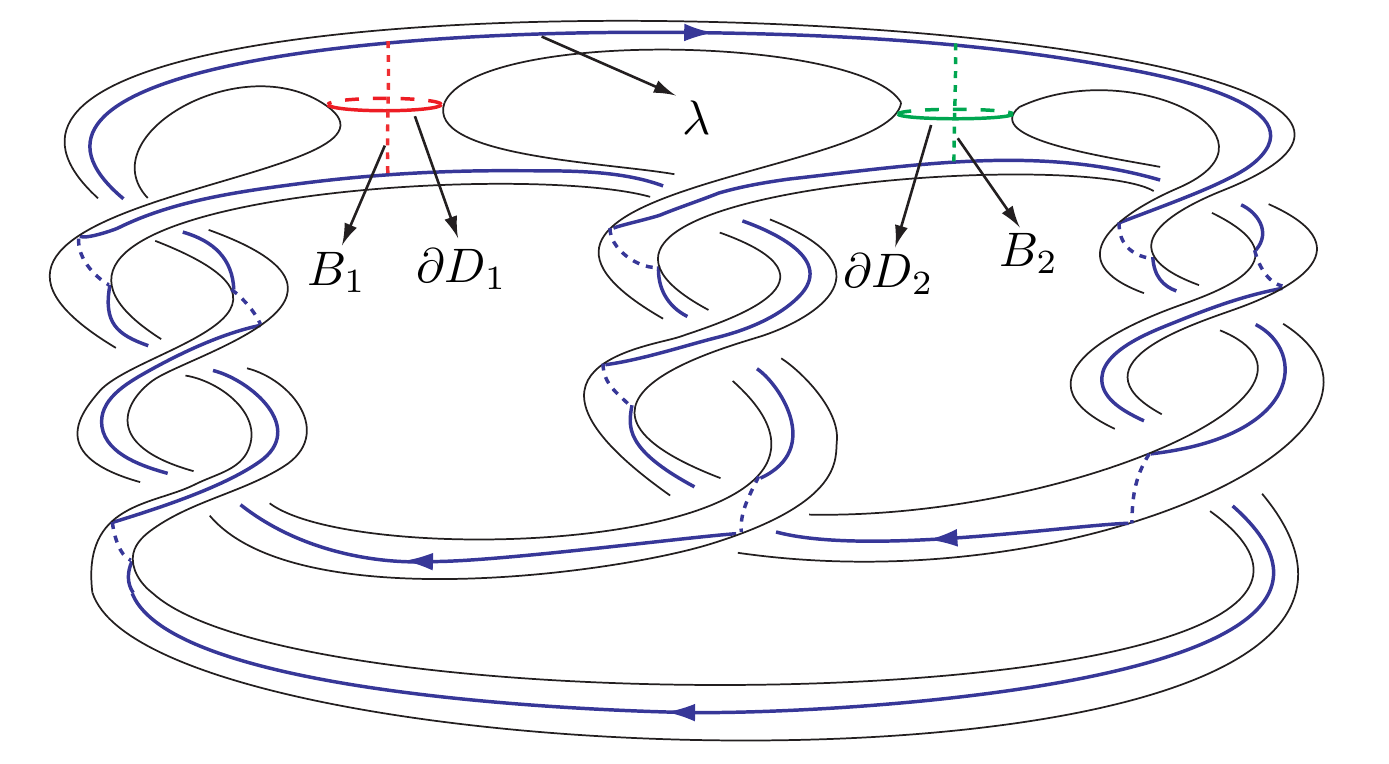}
\end{center}
\caption{\label{333+tunnel} $H$, a regular neighborhood of $K=(3, 3, 3)$ with unknotting tunnels.}
\end{figure}

\begin{figure}
\begin{center}
\includegraphics[width=4.5in]{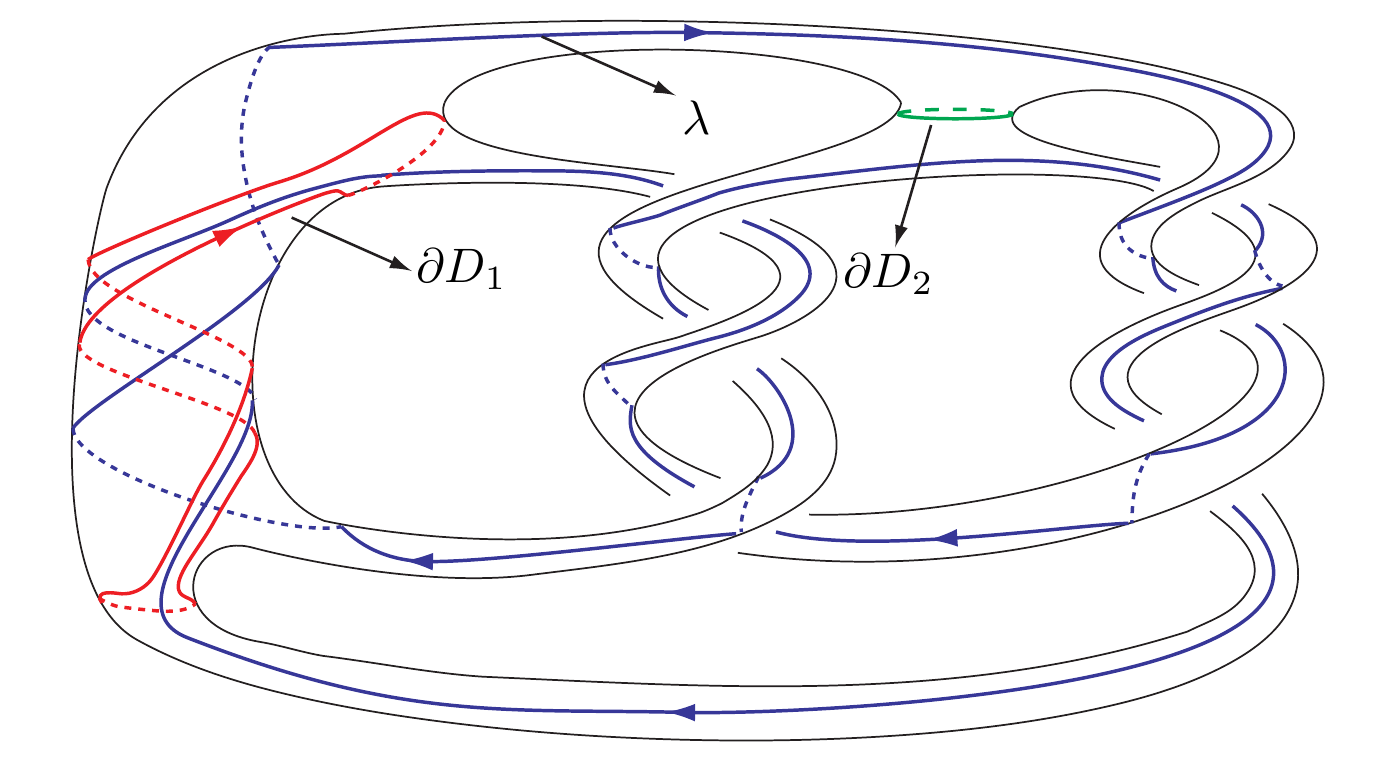}
\end{center}
\caption{\label{3transform1} The deformation of $H$, $\partial D_1$, $\partial D_2$ and $\lambda$ (part
$1$).}
\end{figure}

\begin{figure}
\begin{center}
\includegraphics[width=4.5in]{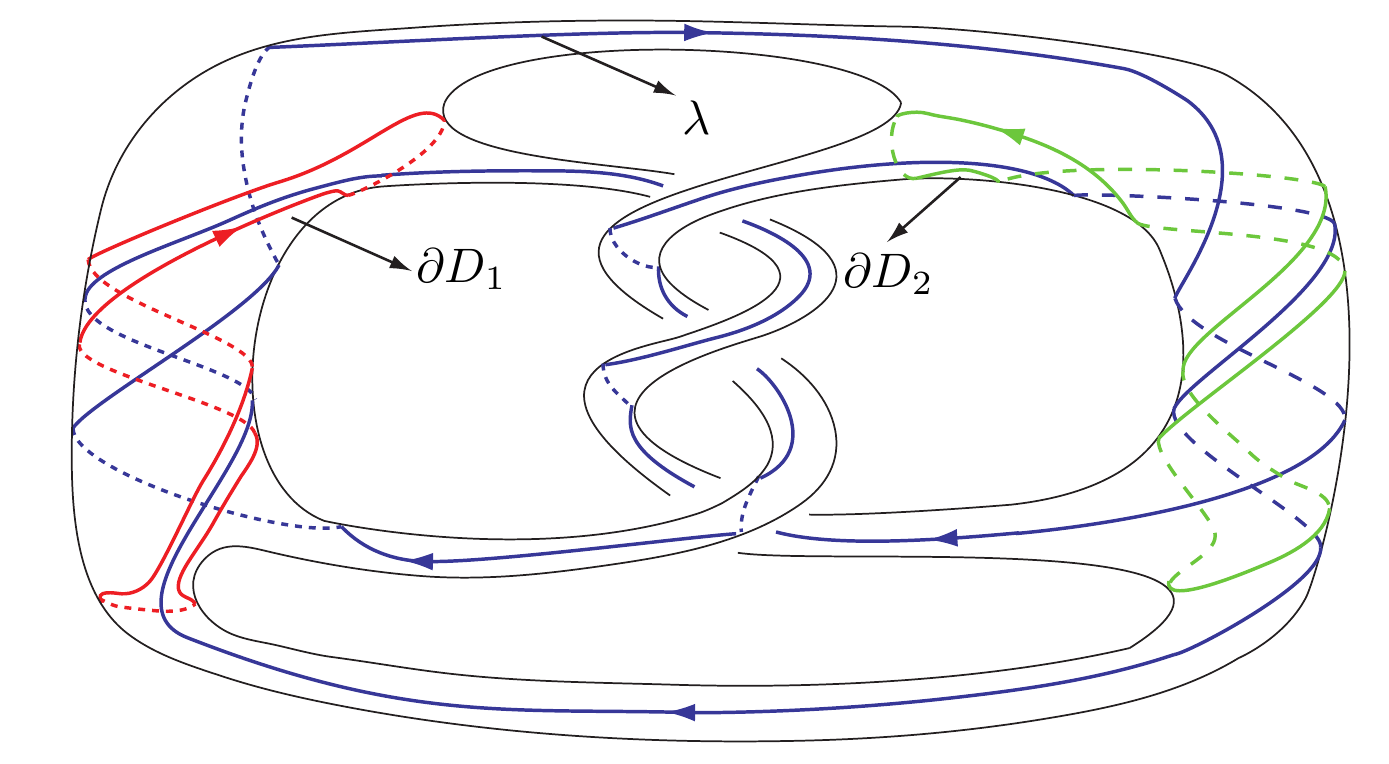}
\end{center}
\caption{\label{3transform2} The deformation of $H$, $\partial D_1$, $\partial D_2$ and $\lambda$ (part
$2$).}
\end{figure}

\begin{figure}
\begin{center}
\includegraphics[width=5in]{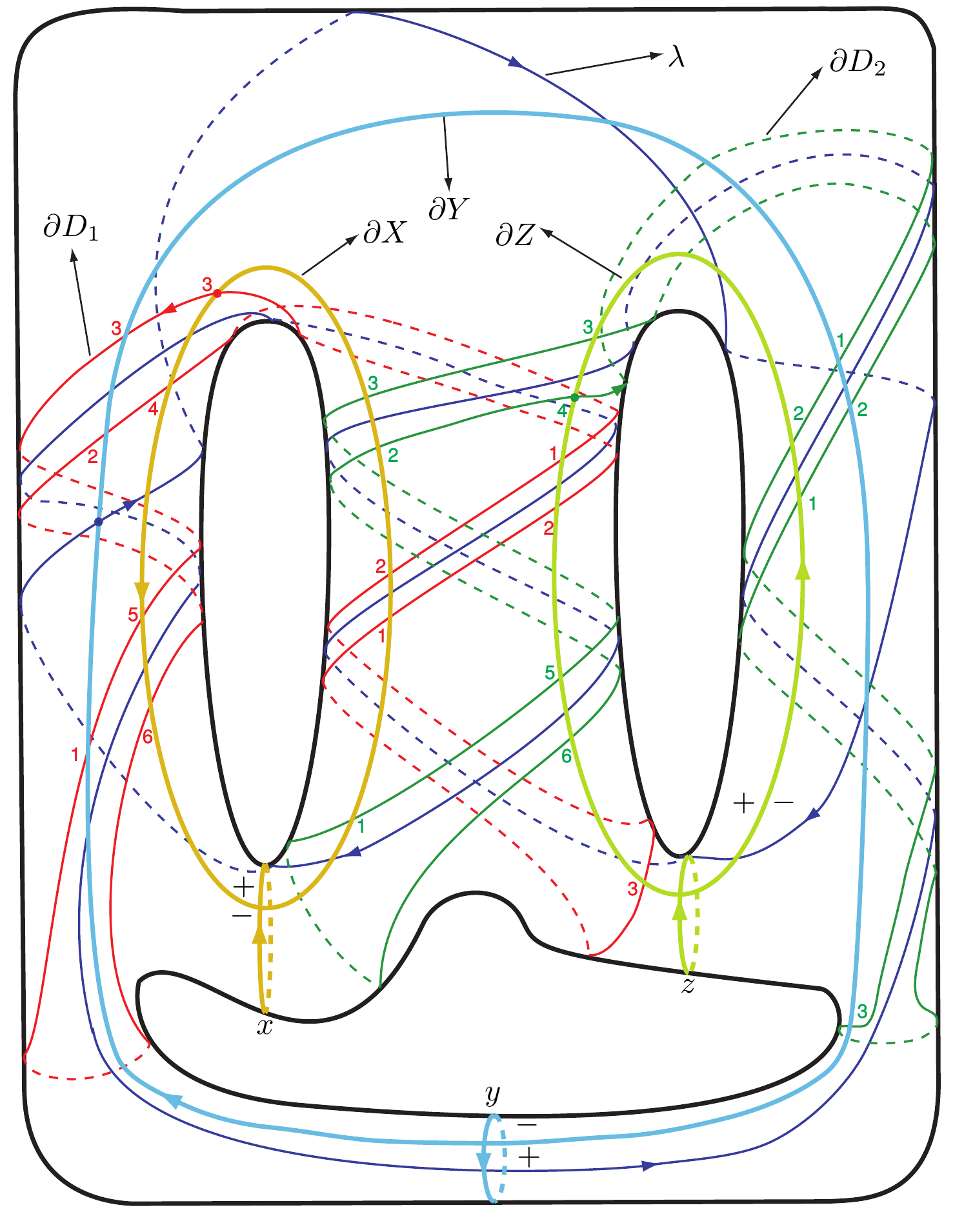}
\end{center}
\caption{\label{3transformed} The deformation of $H$, $\partial D_1$, $\partial D_2$ and $\lambda$ (part
$3$).}
\end{figure}

As showed in \cite{la}, $K$ is a tunnel number two knot.  Figure \ref{333+tunnel} shows two unknotting
tunnels $B_1$ and $B_2$ for $K$ (noticing here we have more options of the tunnels and we always pick the two as shown in Figure \ref{333+tunnel}), and a regular neighborhood $H$ of $K\cup B_1\cup B_2$.
In the
figure, $D_i$ is a meridian disk of $N(B_i)$.
$H$ is a handlebody of genus three. We can deform $H$ such
that its exterior $H'$ is a standard handlebody in $\mathbb{S}^3$. At the same time of the deformation,
we can keep track of the curves $\partial D_1$, $\partial D_2$ and $\lambda$, where $\lambda$ is a
standard longitude. Figures \ref{3transform1}-\ref{3transformed} show the procedure of the deformation:
Figure \ref{3transform1} shows the result after we untangle the three crossings  on the left, Figure
\ref{3transform2} shows the result after  we untangle the three crossings on the right, and Figure
\ref{3transformed} shows the result after we untangle the  three crossings in the middle.

Pick a  disk system $\{X, Y, Z\}$ for $H'$ such that  the boundaries of $X$, $Y$ and $Z$ are as shown in
Figure \ref{3transformed}. Let
 $\{x, y, z\}$ be a generating set for $\pi_1(H')$ dual to disk system, where $x$ is a
simple closed curve in $\partial H'$ which is disjoint from $\partial Y$ and $\partial Z$ and intersects
$\partial X$ exactly once, and $y$ and $z$ are chosen similarly, as indicated in Figure
\ref{3transformed}. Note that in the figure we picked each of $x, y$ and $z$ up to isotopy and  didn't
draw them as loops sharing a common base point.
 We orient $x, y$ and $z$ using the  right-hand rule with  our thumbs pointing to the positive
 direction of $\lambda$. We also assume that $x$ ($y, z$
respectively) travels from the negative side to the positive side of $\partial X$ ($\partial Y, \partial
Z$ respectively). Then we orient $\partial X, \partial Y$ and $\partial Z$ by the right-hand rule such
that our thumbs point to the positive sides. In the figure, we also give indices to the intersection
points between $\{\partial D_1,
\partial D_2\}$ and $\{\partial X, \partial Y, \partial Z\}$ (e.g. $\partial D_1$ (in red color) has six
 intersection points with $\partial X$ and they are labeled  by 1,2,..,6 around $\partial X$
 (in red colo), other intersection points are labeled in similar way).

Following the given directions, we can write out the expressions of $\partial D_1$ and $\partial D_2$ in
terms of $x$, $y$ and $z$. They are
\begin{equation}\label{d1d2}
\parbox{10cm}
{\begin{eqnarray*}
\partial D_1=(x^{-1}y)^2(xy^{-1})(xz^{-1})^2(x^{-1}z),\\
\partial D_2=(zy^{-1})^2(z^{-1}y)(z^{-1}x)^2(zx^{-1}).
\end{eqnarray*}}
\hfill
\end{equation}

So  we get a presentation of the fundamental group $\pi_1(M_K)$ of $M_K$:
\begin{equation}
\pi_1(M_K)=\langle x, y, z : \partial D_1=1, \partial D_2 =1\rangle
\end{equation}
By abelinization we get the homology group of $M_K$
\begin{equation}
H_1(M_K)=\langle y \rangle,
\end{equation}
noticing  that each of $x$ and $z$ is also a generator of $H_1(M_K)$.

By following the given direction, we can also find the expression of the longitude $\lambda$ in terms of
the generators $x$, $y$ and $z$:
\begin{equation}
\lambda = (y^{-1}x)(y^{-1}z)^2(x^{-1}z)(x^{-1}y)^2(z^{-1}y)(z^{-1}x)^2.
\end{equation}

$M_K$ has a Heegaard splitting, $M_K=C \cup_{\partial {H'}} H'$. Where $C$ is a compression body
obtained by attaching two $1$-handles $N(B_1)$ and $N(B_2)$ to the positive boundary
$\partial{M_K}\times [1]$ of $\partial{M_K}\times [0, 1]$. $\partial {H'}$ is the Heegaard  surface (of
genus three)  and $\{\partial D_1, \partial D_2, \partial X, \partial Y,
\partial Z\}$ gives us the Heegaard diagram of this splitting, as
shown in Figure \ref{3transformed}.

Now, let's consider the $3$-fold cyclic cover $M_K^3$ of $M_K$ induced by the homomorphism $h$ from
$\pi_1(M_K)$ to $\mathbb{Z}_3$ factoring  through $H_1(M_K)$:
\begin{equation}\label{homo}
h: \pi_1(M_K) \rightarrow \mathbb{Z}_3; x \mapsto \bar{1}, y \mapsto \bar{1}, z \mapsto \bar{1}.
\end{equation}

\begin{figure}
\begin{center}
\includegraphics[width=4in]{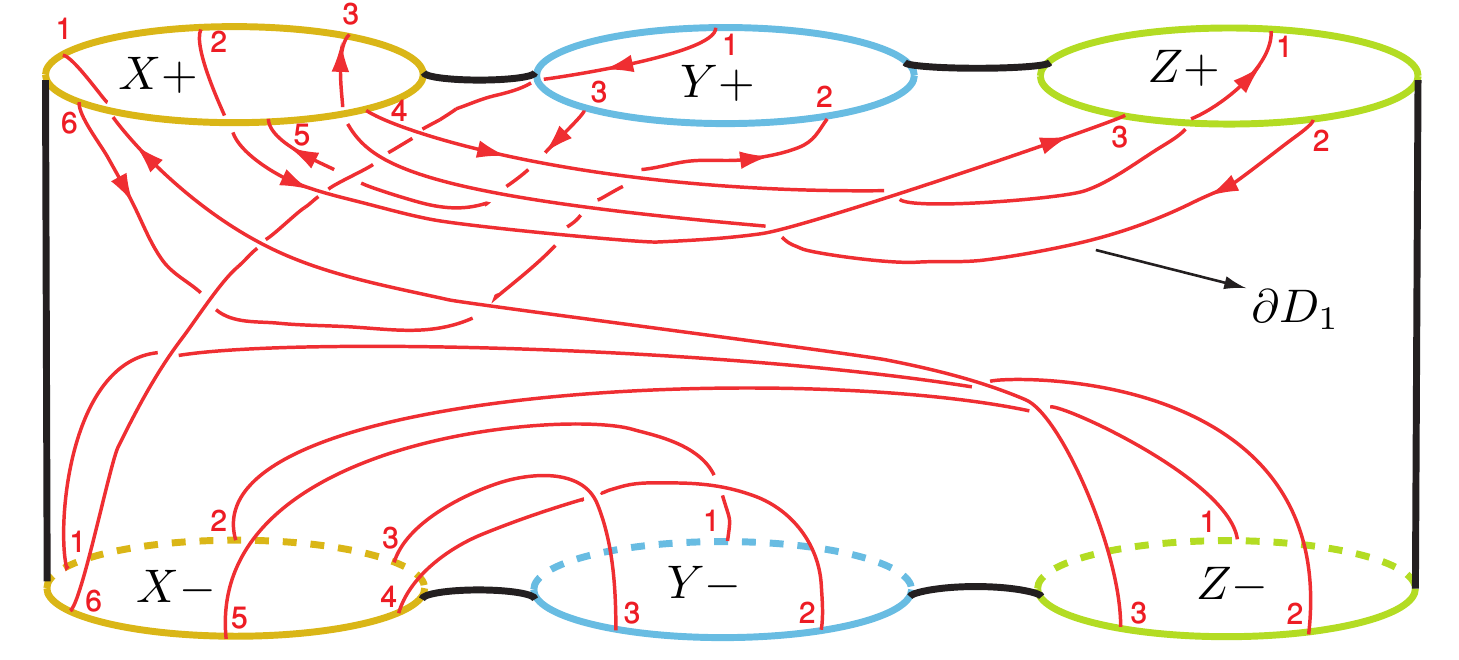}
\end{center}
\caption{\label{333d1} $\partial D_1$ in the resulting surface of cutting $\partial {H'}$ open along
$\{\partial X, \partial Y, \partial Z\}$.}
\end{figure}

\begin{figure}
\begin{center}
\includegraphics[width=5in]{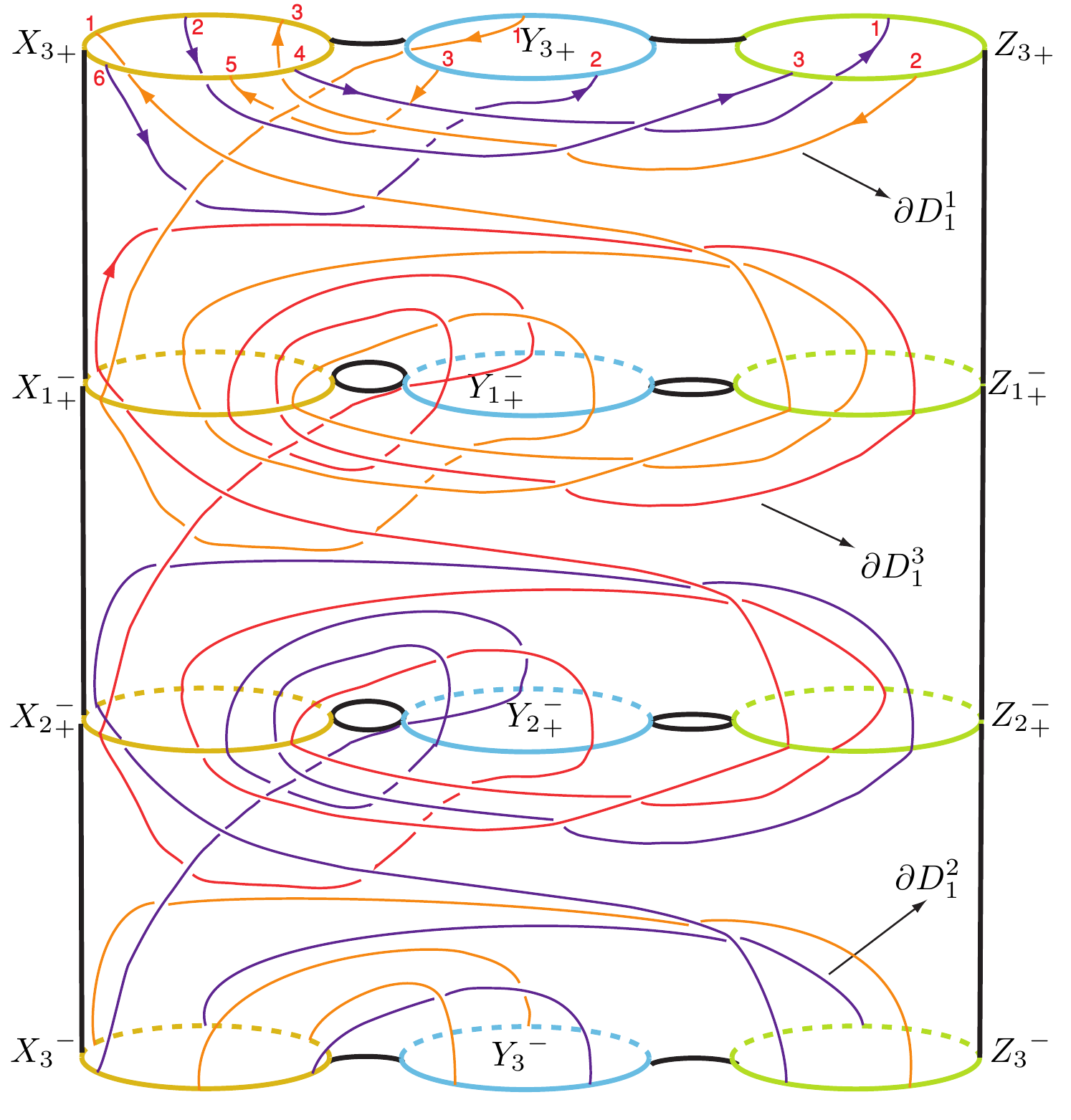}
\end{center}
\caption{\label{333d13} Lifts of $\partial D_1$ in the induced Heegaard surface of $M_K^3$.}
\end{figure}

\begin{figure}
\begin{center}
\includegraphics[width=4in]{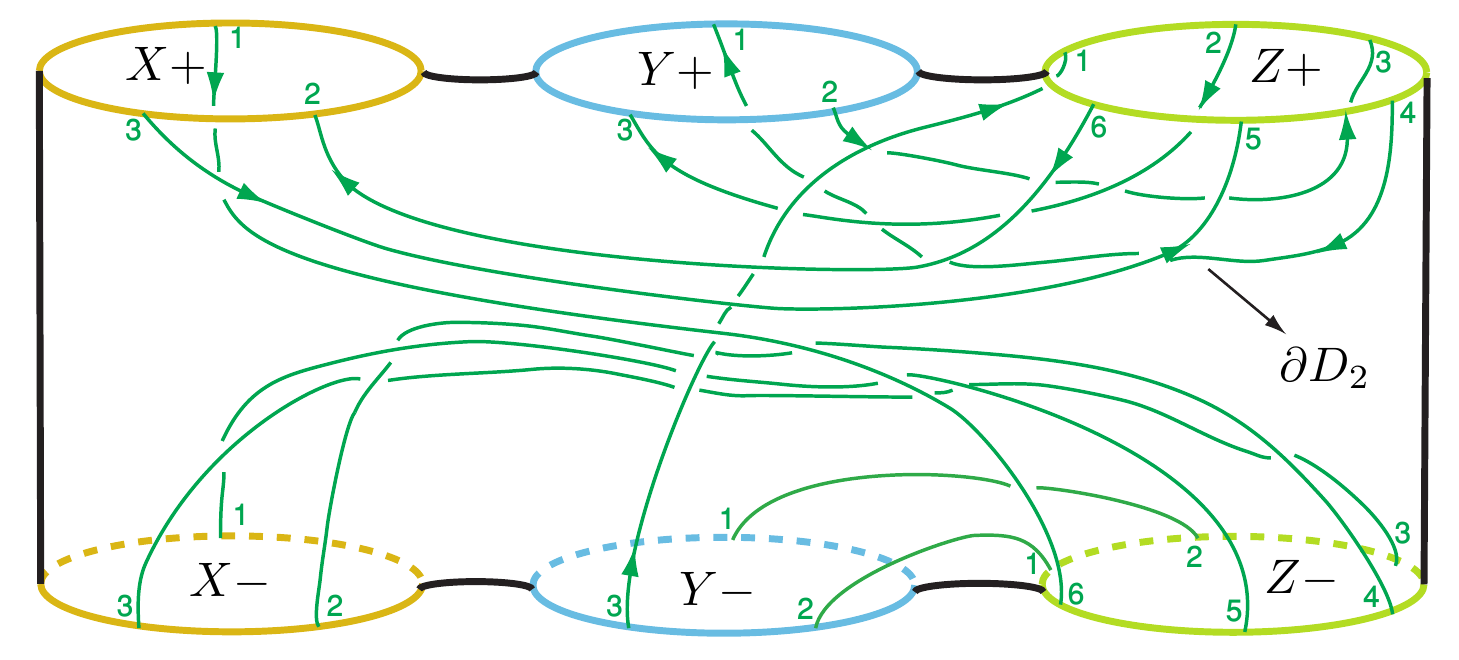}
\end{center}
\caption{\label{333d2} $\partial D_2$ in the resulting surface of cutting $\partial {H'}$ open along
$\{\partial X, \partial Y, \partial Z\}$.}
\end{figure}

\begin{figure}
\begin{center}
\includegraphics[width=5in]{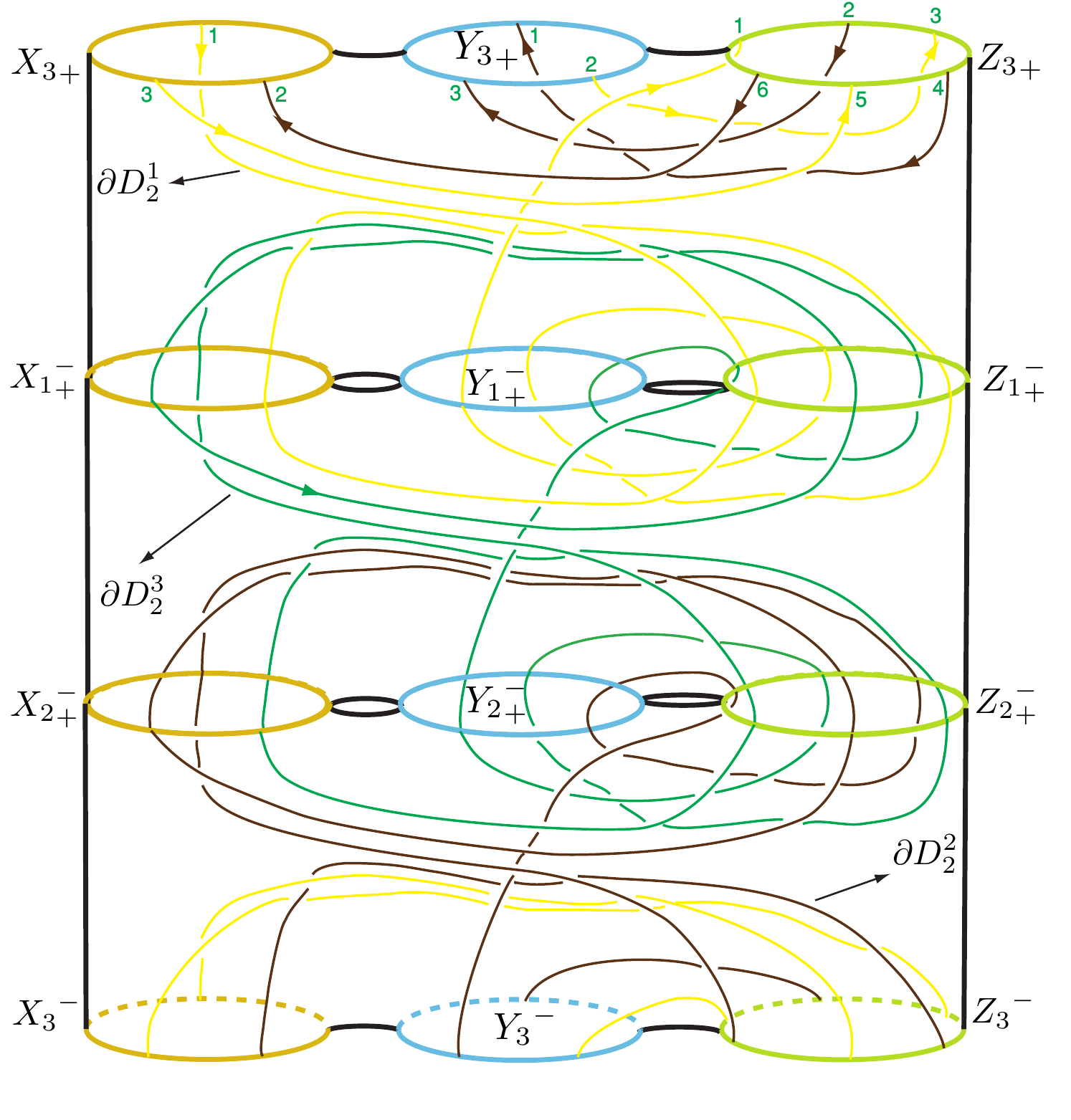}
\end{center}
\caption{\label{333d23} Lifts of $\partial D_2$ in the induced Heegaard surface of $M_K^3$.}
\end{figure}

By cutting $\partial {H'}$ open along $\{\partial X, \partial Y, \partial Z\}$ and pasting $3$ copies of
the resulting surface together cyclicly, we get the induced Heegaard surface of the induced Heegaard splitting of $M_K^3$. We show the
procedure in Figure \ref{333d1}-Figure \ref{333d23}.  Here we should mention that in Figure \ref{333d1}
and Figure \ref{333d2} the curve segments induced from  $\partial  D_i$ are only drawn  schematically.
In reality they
 are embedded on the boundary  surface, but for simplicity,
 we  draw them crossing each other  but keep  their endpoints fixed.
 This simplification will not affect our  proofs because,
later, when we make use of the  Whitehead graphs, we only need information from the  endpoints
 of the curve
segments. Figure \ref{333d13} and Figure \ref{333d23} show us a genus $7$ handlebody $\widetilde{H}$,
which covers  $H'$.  We  take  $\{ X_1, X_2, X_3, Z_1, Z_2, Z_3, Y_3\}$ as a disk system for
$\widetilde{H}$. Each $D_i$  is lifted to three disks $D_i^j$, $j=1,2,3$, whose boundaries are shown in
Figure \ref{333d13} and Figure \ref{333d23}). Let $\widetilde{C}$ be the corresponding cover of $C$,
then $\{ D_1^1, D_1^2, D_1^3, D_2^1, D_2^2, D_2^3\}$ is a disk system for $\widetilde{C}$. $M_K^3$ has
the induced Heegaard splitting $M_K^3=\widetilde{H}\cup_{\partial{\widetilde{H}}}\widetilde{C}$.

\begin{figure}
\begin{center}
\includegraphics[width=4in]{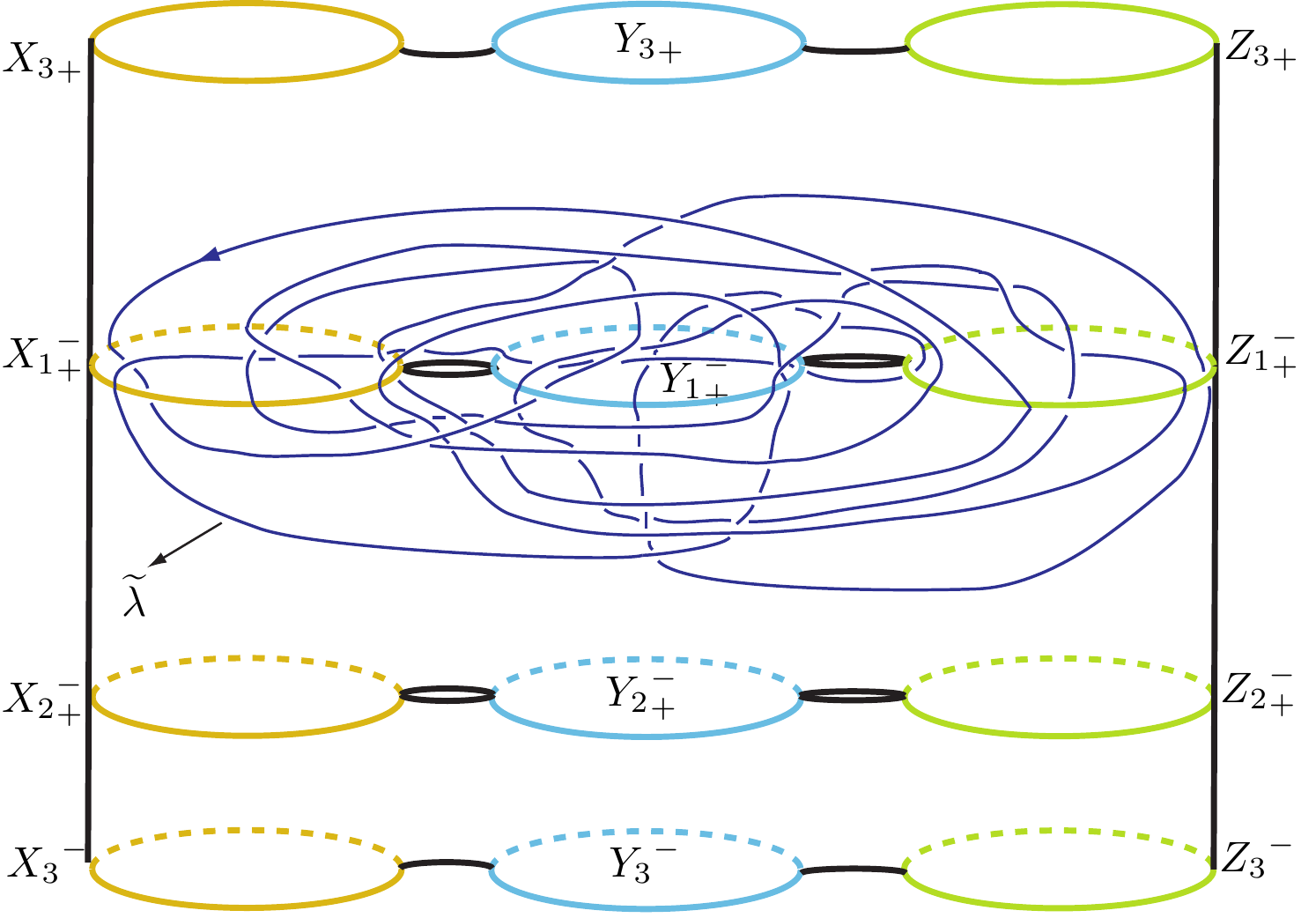}
\end{center}
\caption{\label{3longitude} A lift of $\lambda$ in the induced Heegaard surface of $M_K^3$.}
\end{figure}

From the disk systems we see that the Heegaard splitting of $M_K^3$ is weakly reducible, since $\{
D_1^3, D_2^3\}$ is disjoint form $\{ X_3, Y_3, Z_3\}$. We also notice that the Heegaard splitting is
actually stabilized and thus  reducible. The longitude $\lambda$ is lifted to three copies, we show the
one, $\widetilde{\lambda}$, disjoint from $X_3$, $Y_3$ and $Z_3$ in Figure \ref{3longitude}. Again for
simplicity we did not draw it as embedded on the surface. We note that $\widetilde{\lambda}$ is disjoint
from all $\partial D_i^j$'s.

Now we are going to show  that the closed genus $2$ surface $S$ obtained by compressing the Heegaard
surface $\partial{\widetilde{H}}$ using the disks $D_1^3$, $D_2^3$, $X_3$, $Y_3$ and $Z_3$ is essential
in $M_K^3$. Theorem \ref{theorem 1} will then follow from  Theorem 2.4.3 of \cite{cgls}, namely the surface $S$
remains incompressible in every Dehn filling of $M_K^3$ with slope $m/n$, $(m,n)=1$, $|m|>1$, and since
every such manifold is a cover of the manifold obtained by Dehn filling $M_K$ with slope $3m/n$, $(3m,
n)=1$, $|m|>1$.

\begin{figure}
\begin{center}
\includegraphics[width=4in]{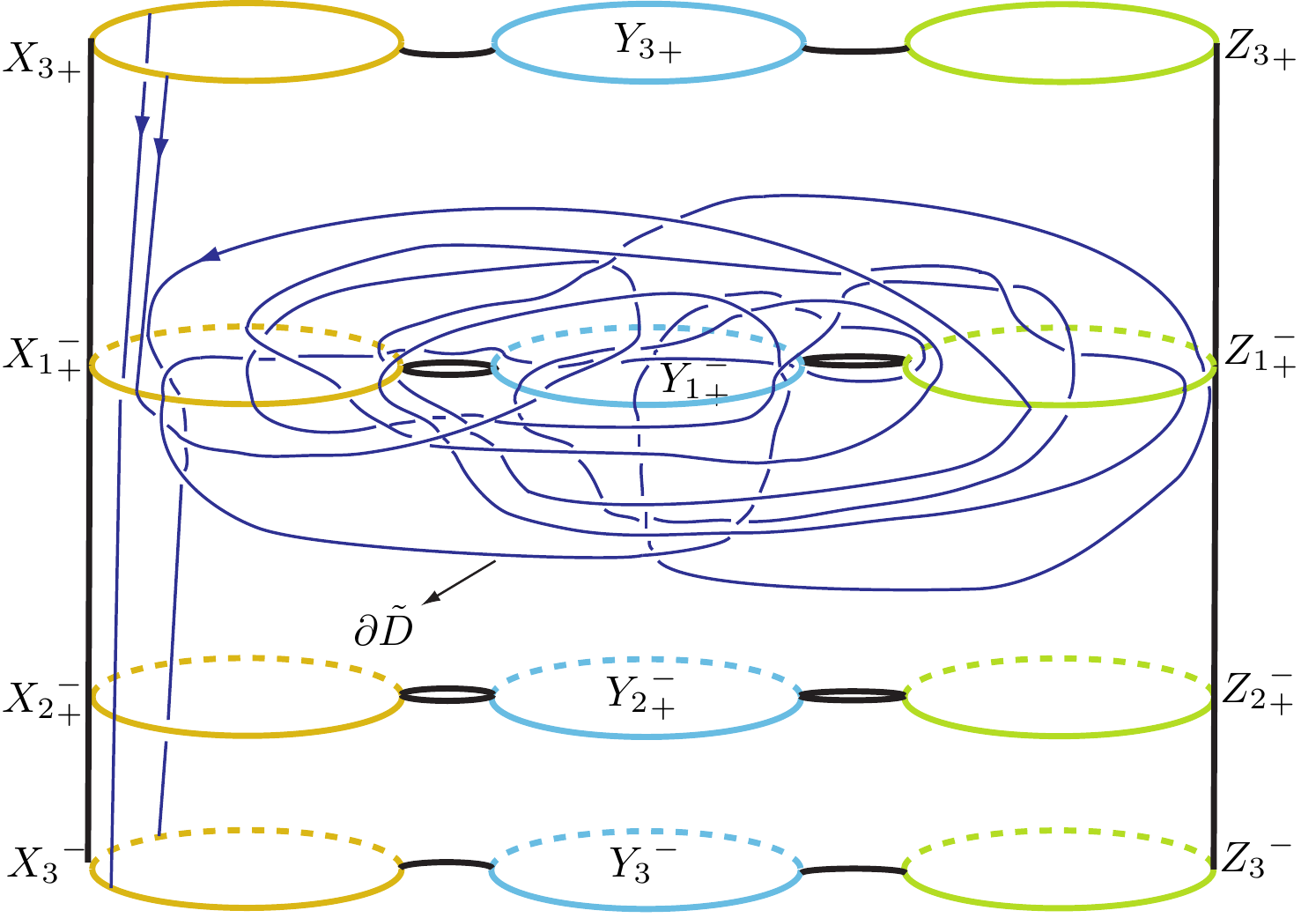}
\end{center}
\caption{\label{3d} $\partial \tilde{D}$ in the derived Heegaard surface of $M_K^3$.}
\end{figure}

It's enough to show that $S$ is incompressible in $M_K^3(2)$, which is the manifold obtained by Dehn
filling $M_K^3$ with slope $2$. Let $C(6)$ be the genus 3 handlebody obtained by Dehn filling $C$ with
slope 6. $M_K^3(2)$ has the induced Heegaard splitting $\widetilde{H}\cup \widetilde{C}(2)$, where
$\widetilde{C}(2)$ is the genus $7$ handlebody covering  $C(6)$. Let $D$ be the meridian disk of the
filling solid torus in $M_K$. By the definition of $D$, we can write out a presentation of $\partial
D$ in terms of $x$, $y$ and $z$:
\begin{equation}
\partial D=\lambda x^6=(y^{-1}x)(y^{-1}z)^2(x^{-1}z)(x^{-1}y)^2(z^{-1}y)(z^{-1}x)^2x^6,
\end{equation}
A sketch of a lift of $\partial D$ in the induced Heeggard surface of $M_K^3$, $\partial\tilde{D}$, is shown in Figure \ref{3d}. Then $\{D_1^1, D_1^2,D_1^3,D_2^1,D_2^2,D_2^3,\tilde{D}\}$ is a disk system of $\tilde{C}(2)$.

\begin{figure}
\begin{center}
\includegraphics[width=4in]{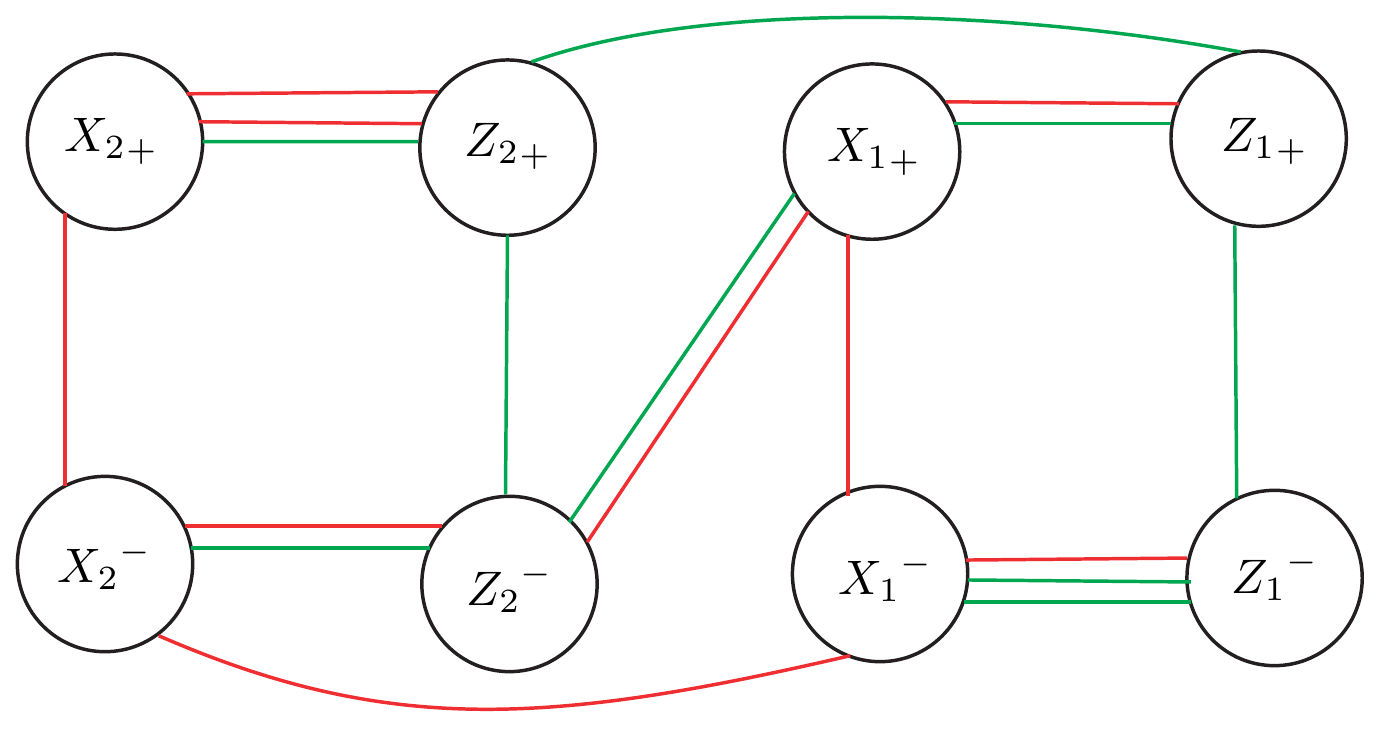}
\end{center}
\caption{\label{3whgraphd1d2} The Whitehead graph of $\{\partial D_1^3, \partial D_2^3\}$ with respect
to $\{X_1, X_2, Z_1, Z_2\}$.}
\end{figure}

Compressing $\widetilde{H}$ along $X_3$, $Y_3$ and $Z_3$, we get a handlebody $\overline{H}$ of genus
$4$, with a disk system $\{X_1, X_2, Z_1, Z_2\}$. $\{\partial D_1^3, \partial D_2^3\}$ is a set of pairwise disjoint
simple closed curves on $\partial \overline{H}$. By following $\partial D_1^3$ and $\partial D_2^3$ on
Figure \ref{333d13} and \ref{333d23}, we can read off the Whitehead graph of $\{\partial D_1^3, \partial
D_2^3\}$ with respect to $\{X_1, X_2, Z_1, Z_2\}$. The graph is shown in Figure \ref{3whgraphd1d2}. The
graph is connected and has no cut vertex. Thus  $\partial{\overline{H}}-(\partial D_1^3\cup
\partial D_2^3)$ is incompressible in $\overline{H}$. If we just look at the Whitehead graph of
$\partial D_1^3$ ($\partial D_2^3$ respectively), i.e., the red part (green part respectively) of the
graph,   there are some valence one vertices in the graph. So by Corollary \ref{lemma 1},   both $\partial
\overline{H}-\partial D_1^3$ and $\partial \overline{H}-\partial D_2^3$ are compressible. Hence by the
multi-handle addition theorem (Theorem \ref{mha}), the manifold $M_1=\overline{H}\cup (D_1^3\times
I)\cup (D_2^3\times I)$ has incompressible boundary.

On the other hand, compressing the handlebody $\widetilde{C}(2)$ along $D_1^3$ and $D_2^3$, we get a
handlebody $\overline{C}$ of genus $5$, with a disk system $\{D_1^1, D_1^2, D_2^1, D_2^2, \tilde{D}\}$.
$\{\partial X_3, \partial Y_3, \partial Z_3\}$ is a set of pairwise disjoint simple closed curves on
$\partial{\overline{C}}$.

\begin{figure}
\begin{center}
\includegraphics[width=4in]{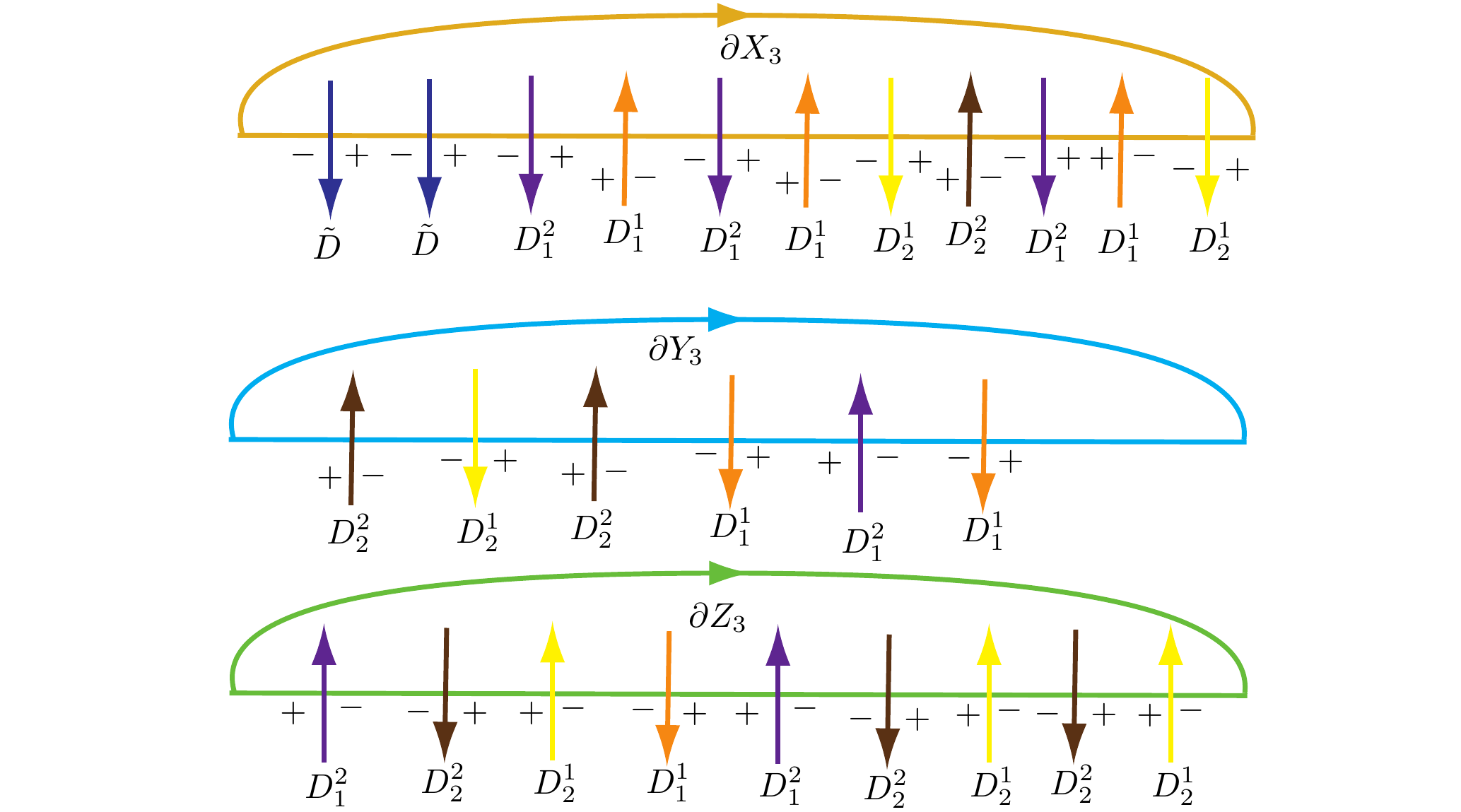}
\end{center}
\caption{\label{xyz} The neighborhoods of $\partial X_3, \partial Y_3$ and $\partial Z_3$.}
\end{figure}

To see the Whitehead graph of $\{\partial X_3, \partial Y_3, \partial Z_3\}$ with respect to $\{D_1^1,
D_1^2, D_2^1, D_2^2, \tilde{D}\}$, we present a neighborhoods of $\partial X_3, \partial Y_3$ and $\partial Z_3$
in Figure \ref{xyz}. In the figure, we mark the positive and negative sides of each of $\{D_1^1, D_1^2,
D_2^1, D_2^2, \tilde{D}\}$ by following the orientations and using right-hand rule, here we always let our
thumbs point to the positive sides.

\begin{figure}
\begin{center}
\includegraphics[width=4in]{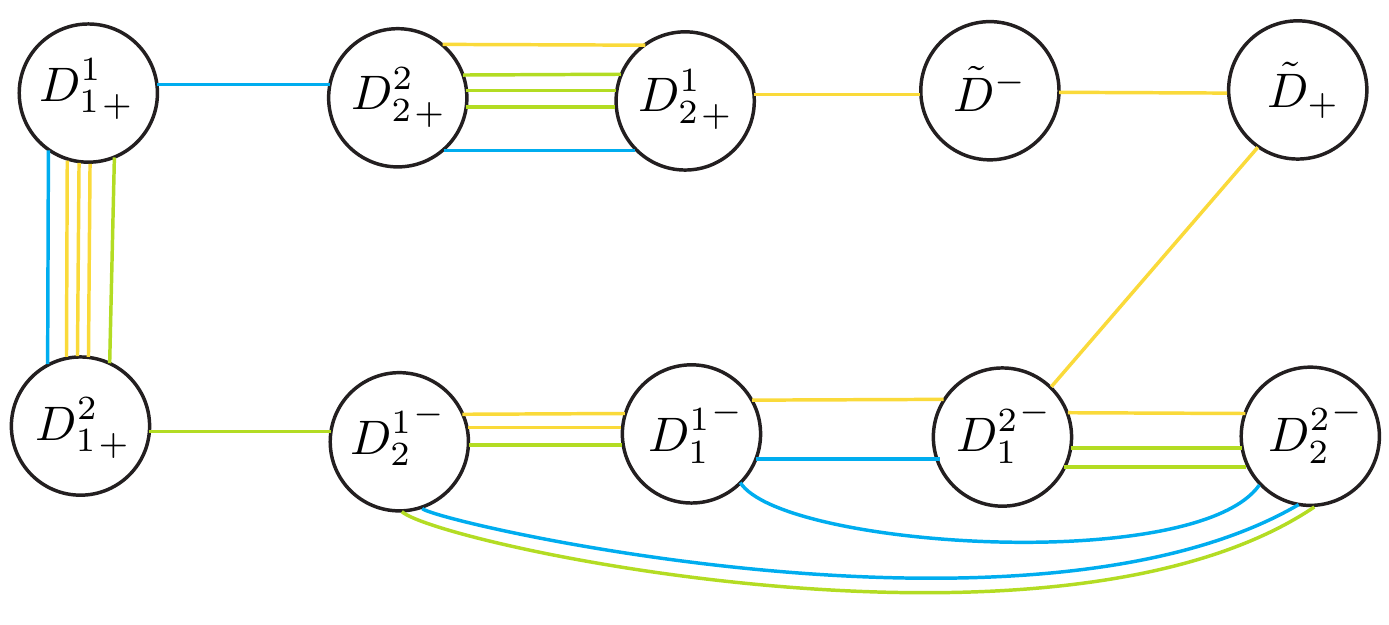}
\end{center}
\caption{\label{3whgraphxyz} The Whitehead graph of $\{\partial X_3, \partial Y_3, \partial Z_3\}$ with
respect to $\{D_1^1, D_1^2, D_2^1, D_2^2, \tilde{D}\}$.}
\end{figure}

By following the oriented  $\partial X_3, \partial Y_3$ and $\partial Z_3$, we can draw the Whitehead
graph of $\{\partial X_3,
\partial Y_3, \partial Z_3\}$ with respect to $\{D_1^1, D_1^2, D_2^1, D_2^2, \tilde{D}\}$, see Figure
\ref{3whgraphxyz}. We can verify   from  this graph that the curve family  $\{\partial X_3,
\partial Y_3, \partial Z_3\}$ satisfies all the conditions of
the multi-handle addition theorem. In fact the graph is connected and has no cut vertex, so
$\partial{\overline{C}}-(\partial X_3\cup \partial Y_3\cup \partial Z_3)$ is incompressible in
$\overline{C}$. If we just look at the graph of $\partial X_3, \partial Y_3$ or $\partial Z_3$, the graph is disconnected, so each of $\partial X_3, \partial Y_3$ and  $\partial Z_3$ does not bind a
free factor $F_4$ of $F_5$. The graph of $\partial
Y_3\cup\partial Z_3$ is disjoint from $D$, so it does not bind $F_5$. The graph of $\partial
X_3\cup\partial Y_3$ (respectively the graph of $\partial X_3\cup\partial Z_3$) has the form as two
subgraphs connected by a path with only two vertices $\{\tilde{D}_-, \tilde{D}_+\}$. By Lemma \ref{lemma 2}, $\{\partial
X_3,
\partial Y_3\}$ (respectively  $\{\partial X_3, \partial Z_3\}$) is
separable in $\partial \overline{C}$, i.e. does not bind $F_5$. So all the conditions of the
multi-handle addition theorem are satisfied. Thus the manifold $M_2=\overline{C}\cup(X_3\times
I)\cup(Y_3\times I)\cup(Z_3\times I)$ has incompressible boundary.

Notice that $S=\partial M_1=\partial M_2$ up to isotopy, so $S$ is incompressible in $M_K^3(2)$. Also
notice that $S$ is contained in $M_K^3$, so $S$ is an essential surface in $M_K^3$. This completes the
proof of Theorem \ref{theorem 1} for the case $K=(3, 3, 3)$.

In general, for a pretzel knot $K=(\pm(2i+1), \pm3, \pm(2j+1))$, the proof is similar.
Up to knot
 equivalence and taking mirror images, we can divide our proof into 3 cases:

(1): $K=(2i+1, 3, 2j+1)$,

(2): $K=(-(2i+1), 3, 2j+1)$,

(3): $K=(-(2i+1), 3, -(2j+1))$.

Let $H$ be a regular neighborhood of $K$ and the unknotting tunnels (which are chosen similarly as we did for the $(3, 3, 3)$-pretzel knot). After
some proper deformation of $H$ in $S^3$, we can make the exterior $H'$ of $H$ in $\mathbb{S}^3$ a
standard handlebody of genus $3$. In the meantime we are still  able to  keep  track of the boundaries
of the meridians of the unknotting tunnels, $\partial D_1$ and $\partial D_2$, and get  their final
appearance in $\partial H'$.

If we pick the disk system $\{X, Y, Z\}$ of $H'$ and the generating set $\{x,y,x\}$ of $\pi_1(H')$ as
before, we can read off a presentation of $\pi_1(M_K)$ from the curves $\partial D_1$ and $\partial
D_2$:
\begin{equation}
\pi_1(M_K)=\langle x, y, z : \partial D_1=1, \partial D_2=1\rangle,
\end{equation}
where $\partial D_1=\begin{cases}
      (x^{-1}y)^{i+1}(xy^{-1})^i(xz^{-1})^2(x^{-1}z)& \text{Case 1}, \\
      (y^{-1}x)^{i}(yx^{-1})^{i+1}(xz^{-1})^2(x^{-1}z)& \text{Case 2}, \\
      (y^{-1}x)^{i}(yx^{-1})^{i+1}(xz^{-1})^2(x^{-1}z)& \text{Case 3};
\end{cases}$

\ \ \ \ $\partial D_2=\begin{cases}
      (zy^{-1})^{j+1}(z^{-1}y)^j(z^{-1}x)^2(zx^{-1})& \text{Case 1}, \\
      (zy^{-1})^{j+1}(z^{-1}y)^{j}(z^{-1}x)^2(zx^{-1})& \text{Case 2}, \\
      (yz^{-1})^{j}(y^{-1}z)^{j+1}(z^{-1}x)^2(zx^{-1})& \text{Case 3}.
 \end{cases}$

By abelinization we get $H_1(M_K)=\mathbb{Z}=\langle y\rangle$, and $x=z=y$.

Similarly we can read off a word expression for an oriented longitude  $\lambda$ as:

$\lambda=\begin{cases}
      (y^{-1}x)^i(y^{-1}z)^{j+1}(x^{-1}z)(x^{-1}y)^{i+1}(z^{-1}y)^j(z^{-1}x)^2& \text{Case 1}, \\
      (y^{-1}x)^i(z^{-1}y)^{j}(z^{-1}x)z^{-1}(yx^{-1})^{i}(zy^{-1})^j(zx^{-1})z& \text{Case 2}, \\
      (x^{-1}z)(y^{-1}x)^{i}(y^{-1}z)^jy^{-1}(xz^{-1})(yx^{-1})^i(yz^{-1})^jy& \text{Case 3}.
\end{cases}$

Similarly, let $M_K^3$ be the $3$-fold cyclic cover of $M_K$ corresponding to the homomorphism
$\pi_1(M_K)\rightarrow H_1(M_K)\rightarrow \mathbb{Z}_3$ with  the induced Heegaard splitting, which is
also weakly reducible and stabilized because $\{X_3, Y_3, Z_3\}$ is disjoint from $\{\partial D_1^3,
\partial D_2^3\}$. We can prove that the closed genus $2$ surface $S$, obtained by compressing the
Heegaard surface of $M_K^3$ using $X_3, Y_3, Z_3, D_1^3$ and $D_2^3$, is essential.

\begin{figure}
\begin{center}
\includegraphics[width=4in]{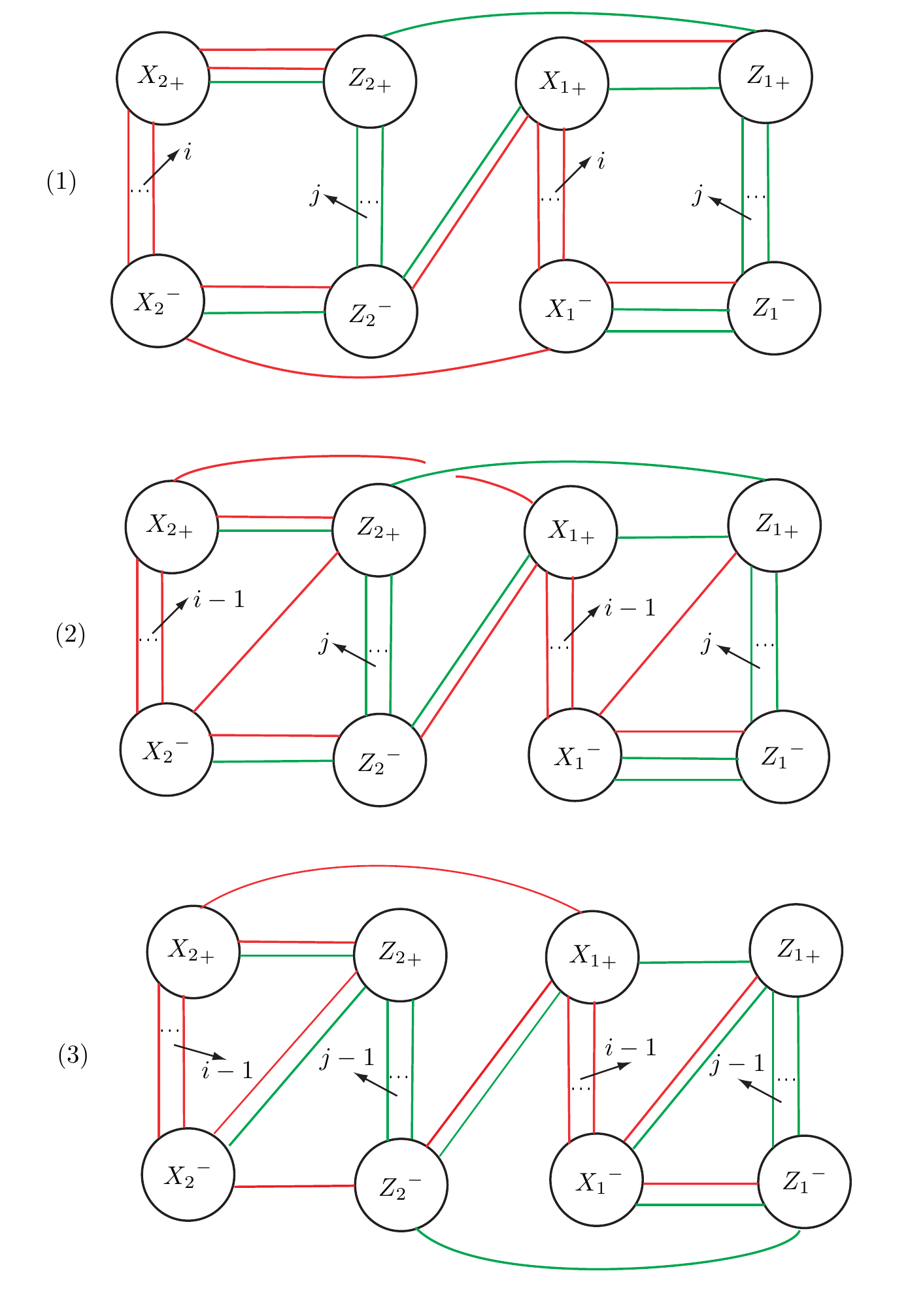}
\end{center}
\caption{\label{whgraphd1d2g} The Whitehead graph of $\{\partial D_1^3, \partial D_2^3\}$ with respect
to $\{X_1, X_2, Z_1, Z_2\}$, (the general case).}
\end{figure}

In fact, the Whitehead graph of $\{\partial D_1^3, \partial D_2^3\}$ with respect to $\{X_1, X_2, Z_1,
Z_2\}$ is as shown in Figure \ref{whgraphd1d2g}, where parts (1)-(3) correspond to Case 1-Case 3
respectively. We can easily check that the graph satisfies all the conditions of the multi-handle
addition theorem.

\begin{figure}
\begin{center}
\includegraphics[width=4in]{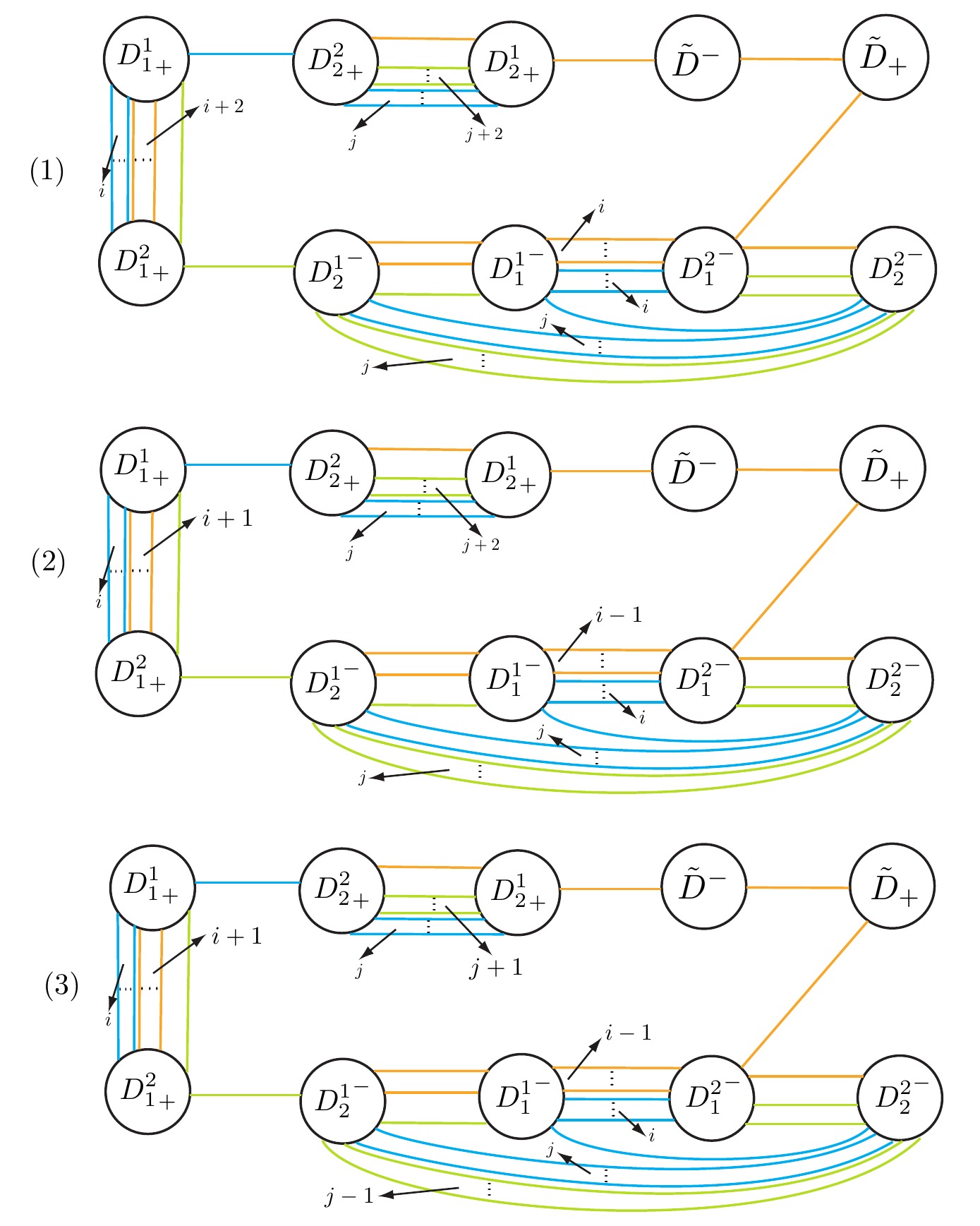}
\end{center}
\caption{\label{whgraphxyzg} The Whitehead graph of $\{\partial X_3, \partial Y_3, \partial Z_3\}$ with
respect to $\{D_1^1, D_1^2, D_2^1, D_2^2, \tilde{D}\}$, (the general case).}
\end{figure}

On the other hand, the Whitehead graph of  $\{\partial X_3, \partial Y_3, \partial Z_3\}$ with respect
to $\{D_1^1, D_1^2, D_2^1, D_2^2, \tilde{D}\}$ is  shown in Figure \ref{whgraphxyzg}, where $\tilde{D}$ is a meridian
disk of the filling torus of the Dehn filling of $M_K^3$ with slope 2, and  parts (1)-(3) correspond to
Case 1-Case 3 respectively. We can also check that the graph satisfies all the conditions of the
multi-handle addition theorem.

This finishes the proof of Theorem \ref{theorem 1}.

\begin{remark}\label{remark 1}
All the lifted Heegaard splittings in Theorem \ref{theorem 1} are stabilized (reducible). We can tell
that from the Heegaard diagrams or from the Whitehead graphs.
\end{remark}

\section{Another application of the method}

In this section, using the same method, we give a new proof of a part of a result of \cite{o}. We show

\begin{prop}\label{prop} If $K=(p_1, p_2, \cdots , p_k)$ is a pretzel knot
with $k\geqslant 4, p_i\geqslant 3$, then $M_K=S^3\setminus K$ contains a closed incompressible
 surface  which remains incompressible in every closed $3$-manifold obtained
 by a non-trivial Dehn filling on $M_K$.
\end{prop}

\begin{remark}\label{remark} Similar method can be used to prove the cases that $p_i\in\mathbb{Z}$ and $|p_i|\geqslant3$.
\end{remark}

According to Formula (\ref{components}) in Section \ref{pre}, a pretzel link $K=(p_1, \cdots, p_k)$ is a
knot if and only if $k$ and all $p_i$'s are odd or exactly one of the $p_i$'s is even. So we may divide
our proof into two cases:

\textbf{Case 1:} $k>4$ odd, all the $p_i$'s are odd.

\textbf{Case 2:} One of the $p_i$'s is even, $k\geqslant 4$.

\textbf{Proof of Case 1:} We will first prove this case for an example, $K=(3, 3, 3, 3, 3)$.

\begin{figure}
\begin{center}
\includegraphics[width=4in]{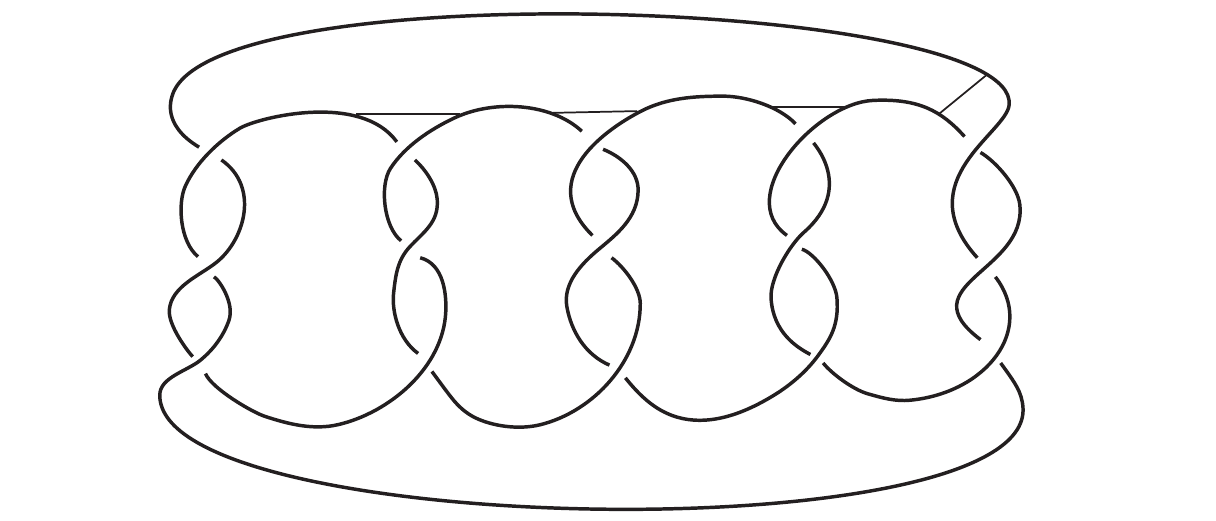}
\end{center}
\caption{\label{33333} The $(3, 3, 3, 3, 3)$-pretzel knot with unknotting tunnels.}
\end{figure}

\begin{figure}
\begin{center}
\includegraphics[width=5.5in]{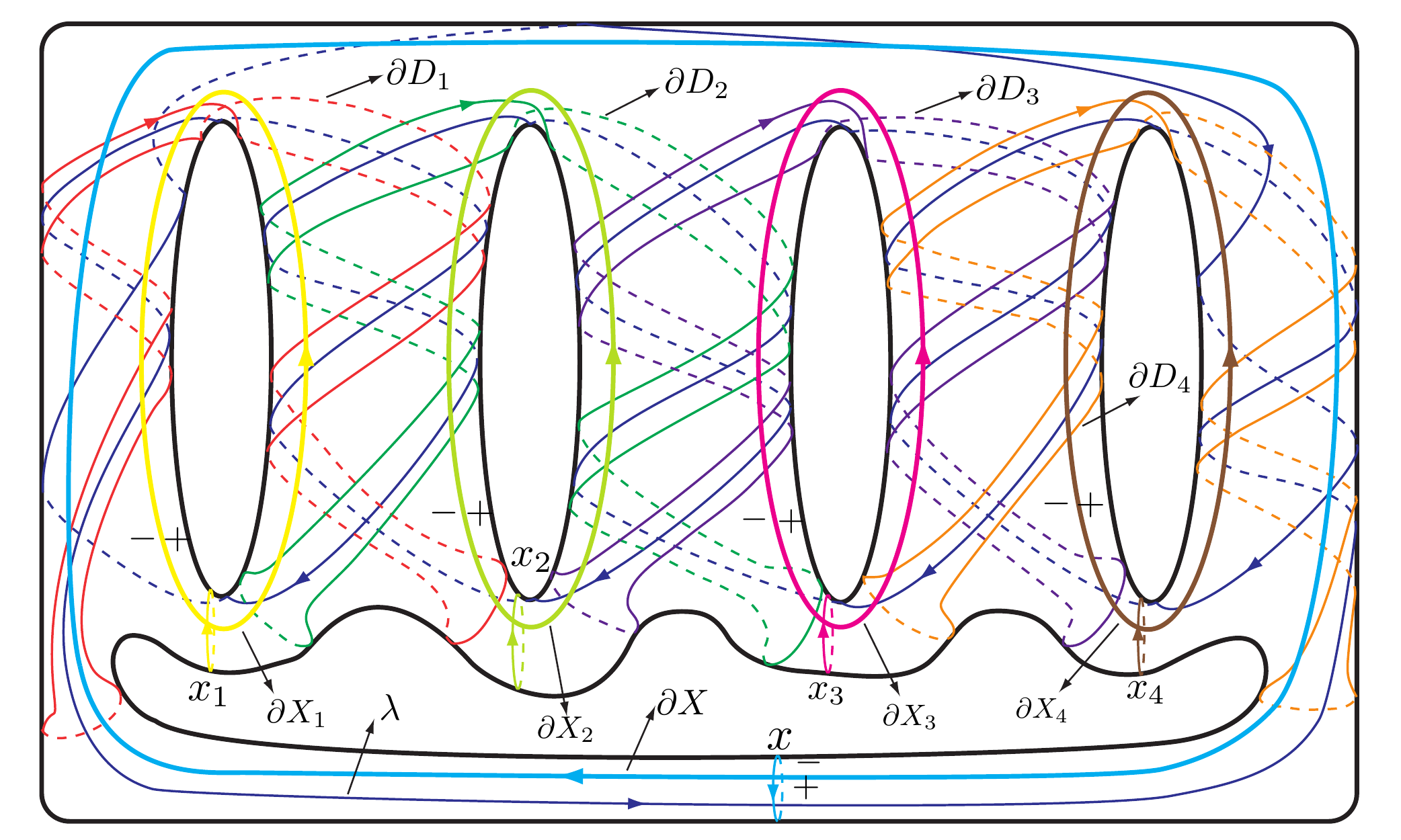}
\end{center}
\caption{\label{33333transformed} The deformation of $H$, $\lambda$, $\partial D_1$, $\partial D_2$,
$\partial D_3$, $\partial D_4$.}
\end{figure}

Figure \ref{33333} shows the $(3, 3, 3, 3, 3)$-pretzel knot $K$ with unknotting tunnels. A regular
neighborhood, $H$, of the union of $K$ and the unknotting tunnels is a genus five handlebody. Let
$D_1,...,D_4$ be meridian disks of the four unknotting tunnels respectively and let $\lambda$ be a
preferred longitude of $K$. We can deform $H$ in $\mathbb{S}^3$ such that its exterior $H'$ is a standard
handlebody in $\mathbb{S}^3$. At the same time we keep tracking the corresponding deformation of the
curves $\partial D_1$, $\partial D_2$, $\partial D_3$, $\partial D_4$ and  $\lambda$. Figure
\ref{33333transformed} shows the final position of these curves on the boundary surface of $H'$. We pick
a disk system $\{X_1, X_2, X_3, X_4, X\}$ of $H'$ and a dual generating set $\{x_1,x_2,x_3,x_4, x\}$ of
$\pi_1(H')$ as shown in Figure \ref{33333transformed}. We orient all the curves in $\partial H'$ by the
same method we used in the proof of Theorem \ref{theorem 1}. Figure \ref{33333transformed} shows us a
Heegaard splitting of the exterior $M_K$ of $K$,  i.e., $M_K=H'\cup C$, where $C$ is a compression body
obtained by attaching four $1$-handles to the positive boundary $\partial{M_K}\times [1]$ of
$\partial{M_K}\times[0, 1]$ with $\{ D_1, D_2, D_3, D_4\}$ as a disk system.

\begin{figure}
\begin{center}
\includegraphics[width=4in]{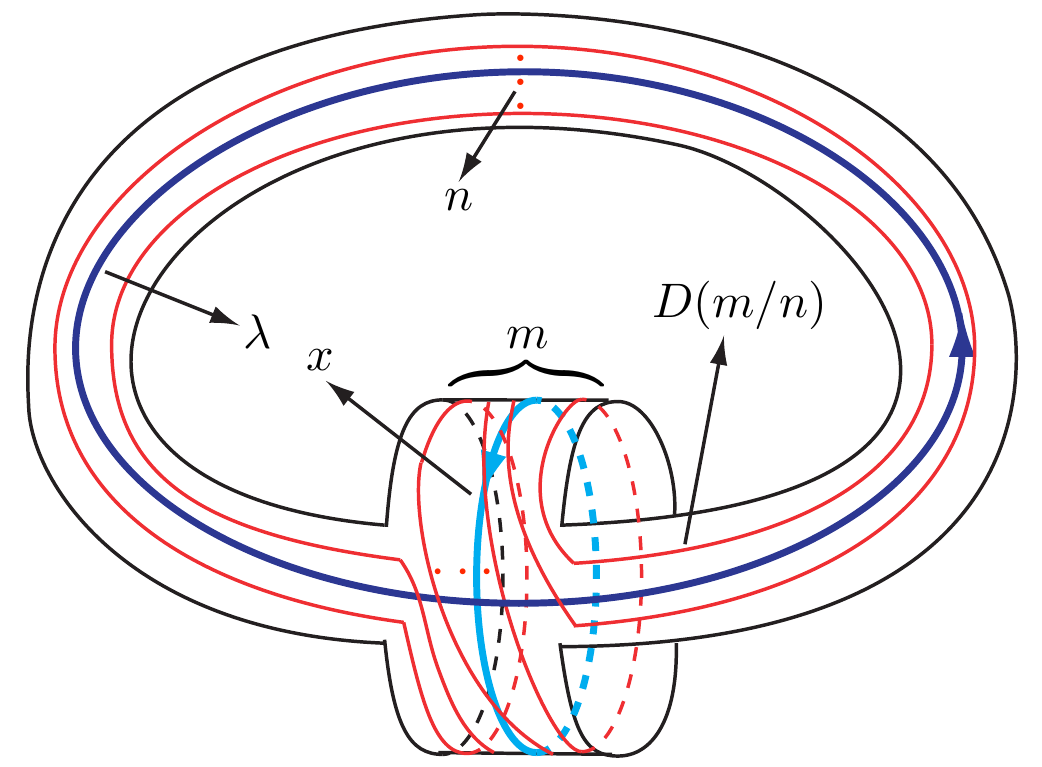}
\end{center}
\caption{\label{Dnm} The filling slope of $C(m/n)$.}
\end{figure}

This Heegaard splitting is weakly reducible, because $\{X\}$ is disjoint from $\{D_2, D_3\}$. We are now
going to show that the genus two surface $S$ obtained by compressing the Heegaard surface $\partial H'$
using $X$, $D_2$ and $D_3$ is incompressible in the manifold $M_K(m/n)$ which is a Dehn filling of $M_K$
with a nontrivial slope $m/n$. The closed manifold $M_K(m/n)$ has the induced  Heegaard splitting
$M_K(m/n)=H'\cup_{\partial H'}C(m/n)$, with $\{ D_1, D_2, D_3, D_4, D(m/n)\}$ as a disk system of
$C(m/n)$ (where $D(m/n)$ is a meridian disk of the Dehn filling torus of $C(m/n)$). Note that $\partial
D(m/n)$ is a simple closed curve on $\partial M_K$, which  can be drawn on a regular neighborhood of
$x\cup\lambda$ as showed in Figure \ref{Dnm}.

\begin{figure}
\begin{center}
\includegraphics[width=4in]{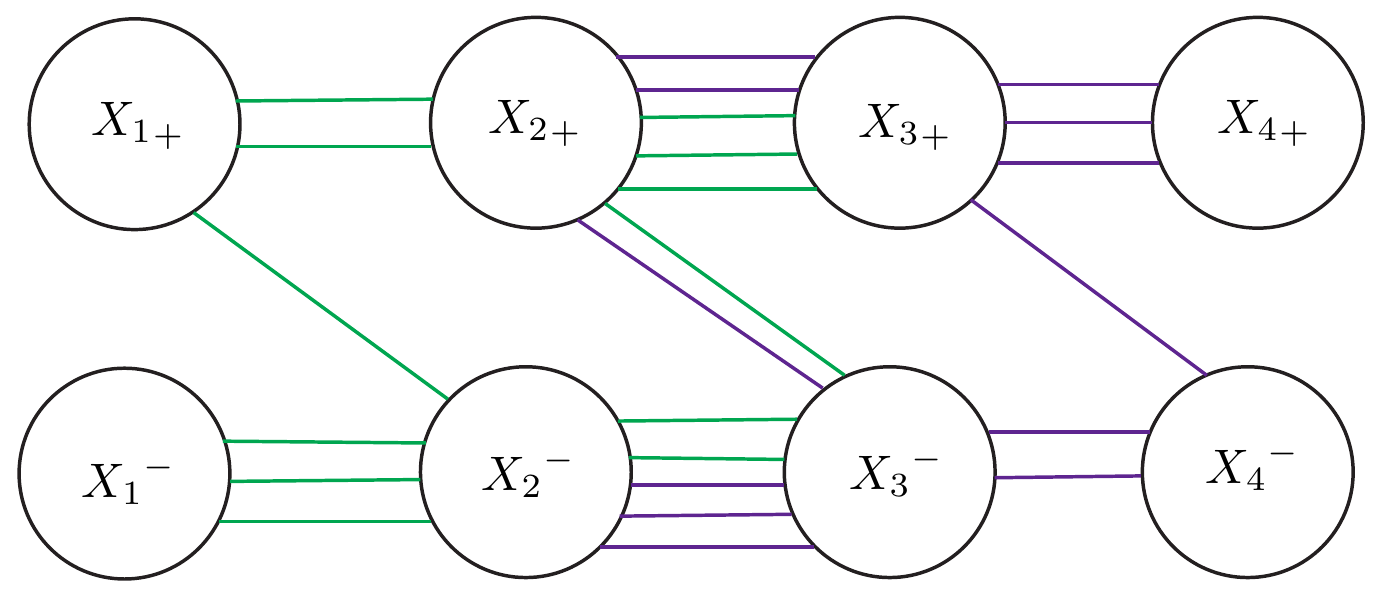}
\end{center}
\caption{\label{whgraphd2d3} The Whitehead graph of $\{\partial D_2, \partial D_3\}$ with respect to
$\{X_1, X_2, X_3, X_4\}$.}
\end{figure}

Compressing $H'$ along $X$, we get a handlebody $\overline{H}$ of genus four, with   disk system $\{
X_1, X_2, X_3, X_4\}$. $\{ \partial D_2, \partial D_3\}$ is a set of simple closed curves on
$\partial{\overline{H}}$. The Whitehead graph of $\{\partial D_2, \partial D_3\}$ with respect to $\{
X_1, X_2, X_3, X_4\}$ is shown in Figure \ref{whgraphd2d3}. The graph is connected and has two cut
vertices ${X_3}_+$ and ${X_2}^-$. Applying the Whitehead automorphisms to these cut vertices (first
$X_3^+$ and then  $X_2^-$), we get the graph shown in  Figure \ref{whgraphd2d32}. 
This graph is
connected and has no cut vertex. So $\partial{\overline{H}}-(\partial D_2\cup \partial D_3)$ is
incompressible in $\overline{H}$. If we just look at the graph of $\partial D_2$
($\partial D_3$ respectively) with respect to $\{ X_1, X_2, X_3, X_4\}$, the graph is disconnected,
which  means $\partial \overline{H}-\partial D_2$ ($\partial \overline{H}-\partial D_3$ respectively) is
compressible. So our graph satisfies all the conditions of the multi-handle addition theorem and thus
the manifold $M_1=\overline{H}\cup (D_2 \times I) \cup (D_3 \times I)$ has incompressible boundary. Note
that $\partial M_1=S$.

\begin{figure}
\begin{center}
\includegraphics[width=4in]{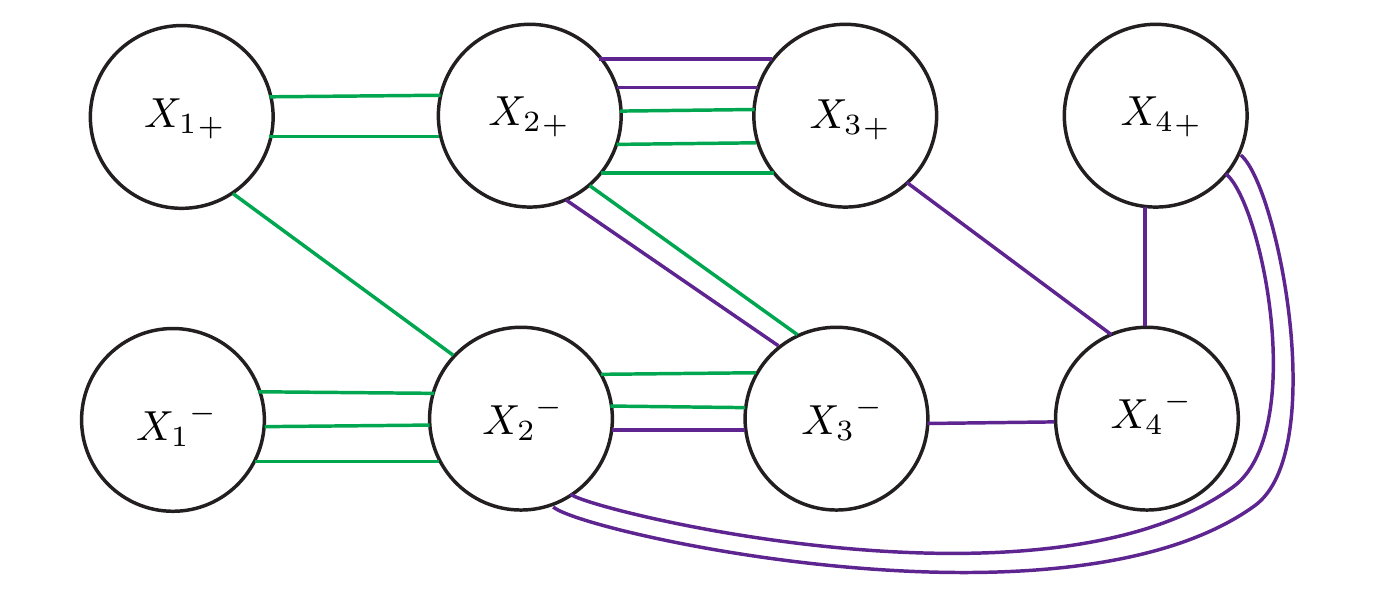}
\end{center}
\caption{\label{whgraphd2d31} The transformed  graph after
 applying the Whitehead algorithm
 to the graph in Figure \ref{whgraphd2d3} at its vertex  ${X_3}_+$.}
\end{figure}

\begin{figure}
\begin{center}
\includegraphics[width=4in]{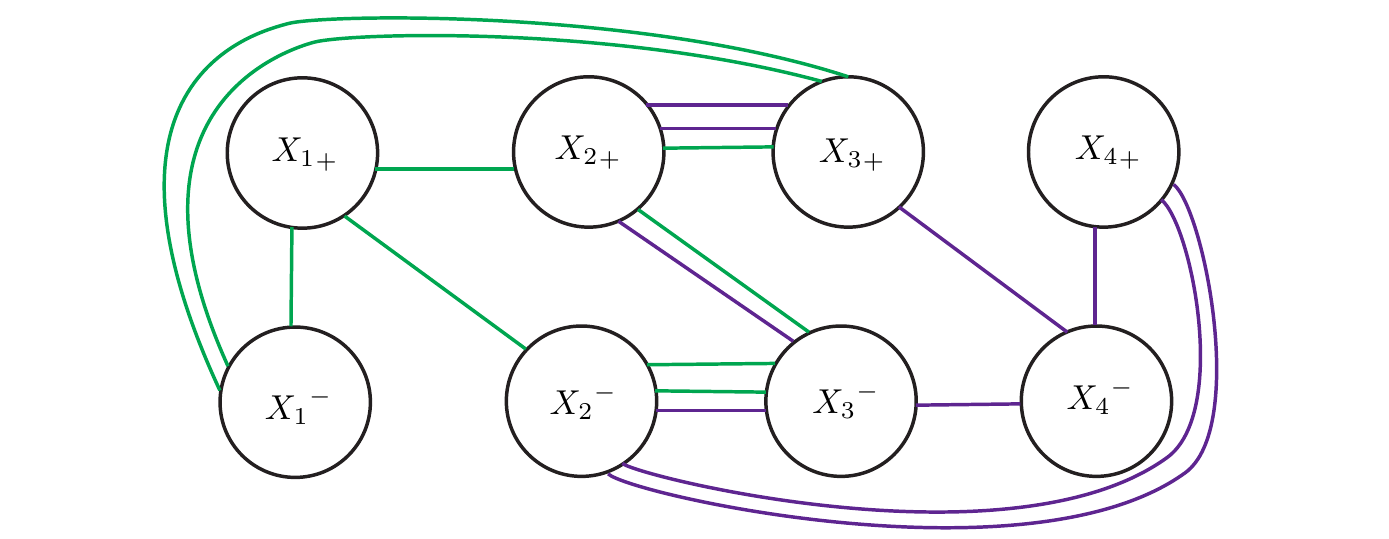}
\end{center}
\caption{\label{whgraphd2d32} The transformed  graph after
 applying the Whitehead algorithm
to the graph in Figure \ref{whgraphd2d31} at its vertex  ${X_2}^-$ .}
\end{figure}

\begin{figure}
\begin{center}
\includegraphics[width=2.5in]{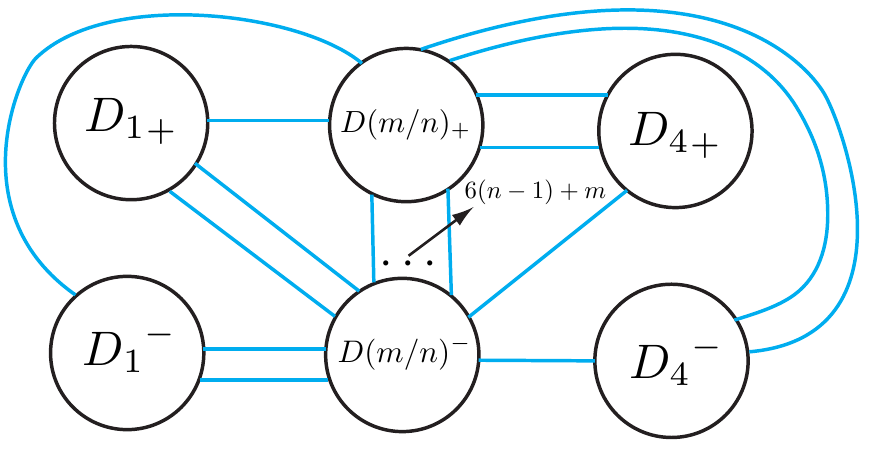}
\end{center}
\caption{\label{whgraphx} The Whitehead graph of $\partial X$ with respect to $\{D_1, D_4, D(m/n)\}$.}
\end{figure}

On the other hand, compress $C(m/n)$ along $D_2$ and $D_3$, we get a handlebody $\overline{C}$ of genus
three with $\{ D_1, D_4, D(m/n)\}$ as a   disk system (for any $m/n\ne 1/0$). $\partial X$ is a simple
closed curve on $\partial{\overline{C}}$. The Whitehead graph of $\partial X$ with respect of $\{ D_1,
D_4, D(m/n)\}$ is shown in Figure \ref{whgraphx}. The graph is connected and has no cut vertex. So
$\partial{\overline{C}}-\partial X$ is incompressible in $\overline{C}$. Hence  the manifold
$M_2=\overline{C}\cup (X\times I)$ has incompressible boundary. Note that  $\partial M_2=S$.

Now we have shown $S=\partial M_1=\partial M_2$ is incompressible in $M_K(m/n)$. Notice that $S$ is
contained in $M_K$, so $S$ is also an essential surface in $M_K$. We finish the proof of Case 1 for the
example $K=(3, 3, 3, 3, 3)$.

In general, for a pretzel knot $K=(p_1, p_2, \cdots, p_k)$ as in Case 1, the proof is similar.

A regular neighborhood, $H$, of the union of $K$ and $k-1$ unknotting tunnels (similarly chosen as we did for the example) is a
genus $k$ handlebody. We can deform $H$ so that its exterior $H'$ is a standard handlebody in
$\mathbb{S}^3$, and at the same time we may track the boundaries of the meridian disks of the unknotting
tunnels, $\partial D_1$, $\cdots$, $\partial D_{k-1}$ and a preferred longitude $\lambda$. We pick a
 disk system of $H'$ and a generating set of $\pi_1(H')$ in a similar way.
 The complement of $K$, $M_K$, has a Heegaard
splitting, $M_K=H'\cup C$, where $C$ is a compression body obtained by attaching $k-1$ $1$-handles to
the positive boundary $\partial{M_K}\times [1]$ of $\partial{M_K}\times[0, 1]$. $\{ X, X_1, \cdots,
X_{k-1}\}$ is a   disk system of $H'$, and $\{ D_1, \cdots, D_{k-1}\}$ is a   disk system of $C$.

Since $\{X\}$ is disjoint from $\{D_2, \cdots, D_{k-2}\}$, this splitting is weakly reducible. We can
show that the genus two surface $S$ obtained by compressing the Heegaard surface $\partial H'$ using
$X$, $D_2, \cdots, D_{k-2}$ is essential in the manifold $M_K(m/n)$ for every  $m/n\ne 1$. $M_K(m/n)$
has a Heegaard splitting $M_K(m/n)=H'\cup_{\partial H'}C(m/n)$. $C(m/n)$ is a genus $k$ handlebody with
a meridian disk system $\{ D_1, \cdots, D_{k-1}, D(m/n)\}$, where $D(m/n)$ is a meridian disk of the
Dehn filling torus of $C(m/n)$.

\begin{figure}
\begin{center}
\includegraphics[width=5.5in]{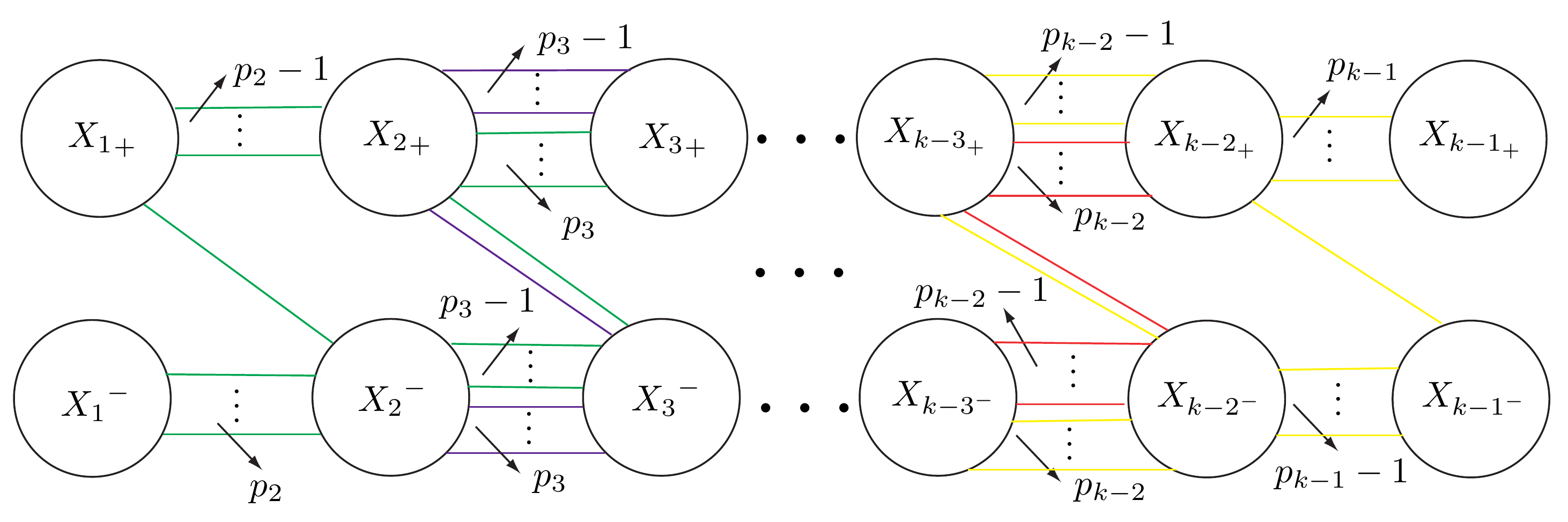}
\end{center}
\caption{\label{whgraphd2dk-2} The Whitehead graph of $\{\partial D_2, \cdots, \partial D_{k-2}\}$ with
respect to $\{X_1, \cdots, X_{k-1}\}$.}
\end{figure}

Compressing $H'$ along $X$, we get a handlebody $\overline{H}$ of genus $k-1$, with  disk system $\{
X_1, \cdots, X_{k-1}\}$. $\{ \partial D_2, \cdots, \partial D_{k-2}\}$ is a set of simple closed curves
on $\partial{\overline{H}}$. The Whitehead graph of $\{\partial D_2, \cdots, \partial D_{k-2}\}$ with
respect to $\{ X_1, \cdots, X_{k-1}\}$ is shown in Figure \ref{whgraphd2dk-2}. We can check that this
graph satisfies all the conditions of the multi-handle addition theorem (after applying some Whitehead automorphisms). So the manifold
$M_1=\overline{H}\cup (D_2\times I)\cup\cdots\cup(D_{k-2}\times I)$ has incompressible boundary $S$.

\begin{figure}
\begin{center}
\includegraphics[width=5in]{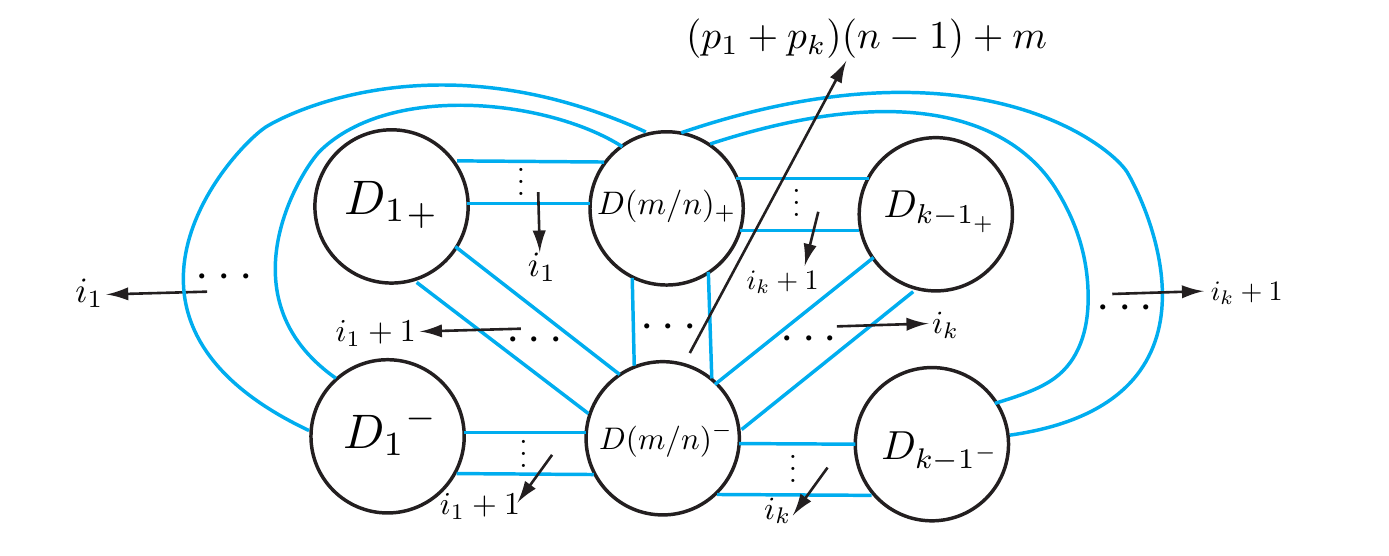}
\end{center}
\caption{\label{whgraphxg} The Whitehead graph of $\partial X$ with respect $\{D_1, D_{k-1}, D(m/n)\}$.}
\end{figure}

On the other hand, compress $C(m/n)$ along $D_2, \cdots, D_{k-2}$, we get a handlebody $\overline{C}$ of
genus three with $\{ D_1, D_{k-1}, D(m/n)\}$ as a  disk system. $\partial X$ is a simple closed curve on
$\partial{\overline{C}}$. The Whitehead graph of $\partial X$ with respect of $\{ D_1, D_{k-1},
D(m/n)\}$ is shown in Figure \ref{whgraphxg}. In the figure, $p_1=2i_1+1$ and $p_k=2i_k+1$. The graph is
connected and has no cut vertex. So $\partial{\overline{C}}-\partial X$ is incompressible in
$\overline{C}$. Thus the manifold $M_2=\overline{C}\cup (X\times I)$ has incompressible boundary.

So  $S=\partial M_1=\partial M_2$ is incompressible in $M_K(m/n)$. Notice that $S$ is contained in
$M_K$, so $S$ is also an essential surface in $M_K$. We finish the proof of Case 1

\begin{figure}
\begin{center}
\includegraphics[width=4in]{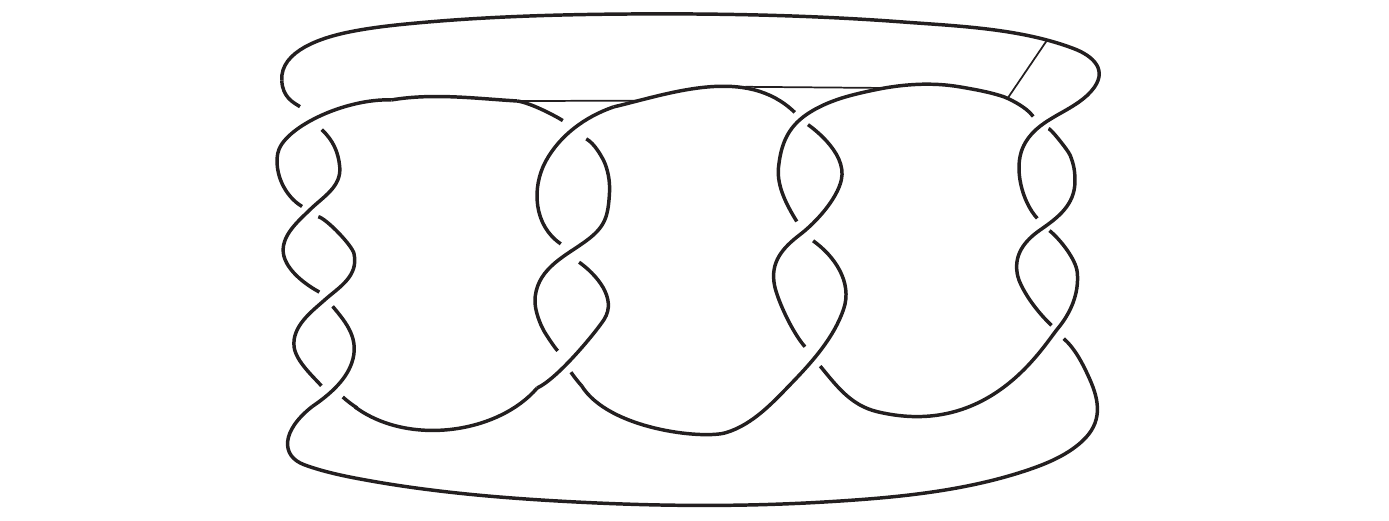}
\end{center}
\caption{\label{3433} The $(4, 3, 3, 3)$-pretzel knot with unknotting tunnels.}
\end{figure}

\textbf{Proof of Case 2:} As before, we first prove this case for an example $K=(4, 3, 3, 3)$. Note
that, up to knot equivalence, we may assume that the left most tangle has  even number of twists. Figure
\ref{3433} shows us the $(4, 3, 3, 3)$-pretzel knot with unknotting tunnels. A regular neighborhood,
$H$, is a handle body of genus four.

\begin{figure}
\begin{center}
\includegraphics[width=5.5in]{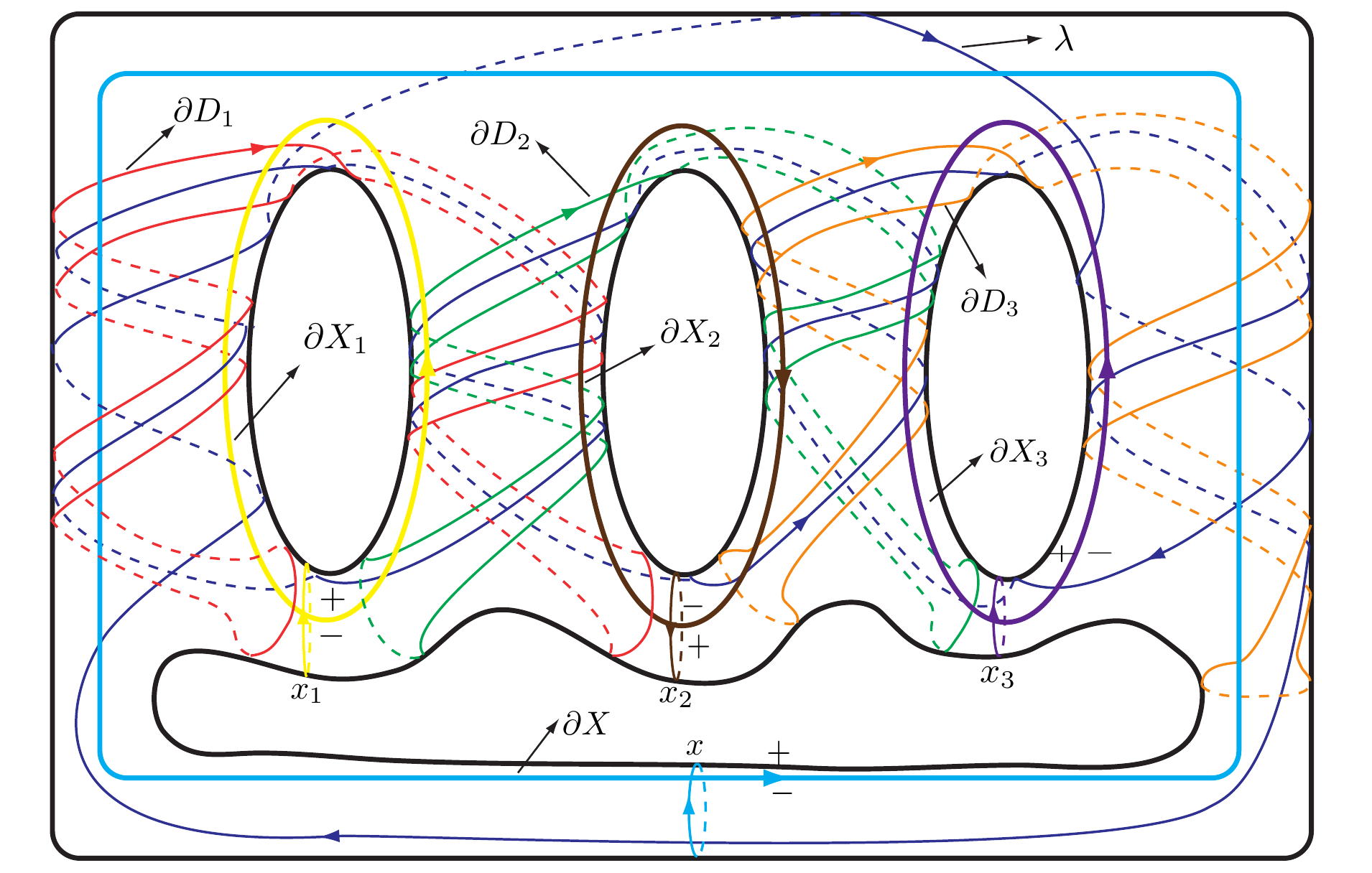}
\end{center}
\caption{\label{3433transformed} The deformation of $H$, $\lambda$, $\partial D_1$, $\partial D_2$,
$\partial D_3$.}
\end{figure}

As before, we deform $H$ so that its exterior looks like a standard handdlebody in $\mathbb{S}^3$. Figure
\ref{3433transformed} shows us the boundary surface of $H$, and $\partial D_1$, $\partial D_2$,
$\partial D_3$, $\lambda$, after the deformation.

The knot complement $M_K=\mathbb{S}^3\setminus K$ has a Heegaard splitting $M_K=H'\cup C$, where $H'$ is
the exterior of $H$ and $C$ is a compression body obtained by attaching three $1$-handles to the
positive boundary $\partial{M_K}\times[1]$. $\{X, X_1, X_2, X_3\}$ is a disk system of $H'$, and
    $\{\partial D_1, \partial D_2, \partial D_3\}$ is a disk system of $C$.

This Heegaard splitting is weakly reducible, because $\{X\}$ is disjoint from $\{D_2\}$. By the similar
argument, one  can show that the genus two surface $S$ obtained by compressing $\partial H'$ using $X$
and $D_2$ is essential in $M_K$. In fact  $S$ remains incompressible in the manifold $M_K(m/n)$ for any
$m/n\ne 1/0$.

\begin{figure}
\begin{center}
\includegraphics[width=2.5in]{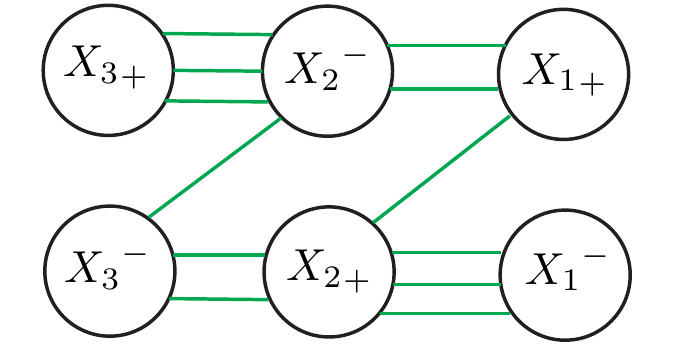}
\end{center}
\caption{\label{whgraphd2} The Whitehead graph of $\partial D_2$ with respect to $\{X_1, X_2, X_3\}$.}
\end{figure}

\begin{figure}
\begin{center}
\includegraphics[width=2.5in]{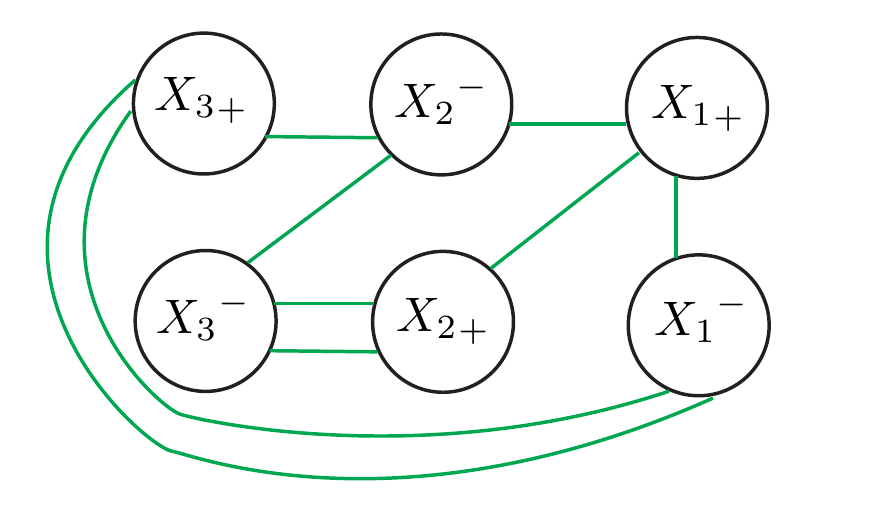}
\end{center}
\caption{\label{whgraphd21} The resulting  graph after  applying the Whitehead automorphism to the graph
in Figure \ref{whgraphd2} at its vertex ${X_2}_+$.}
\end{figure}

Let $\overline{H}$ be the genus three handlebody obtained by compressing $H'$ along $X$. It has a
  disk system $\{X_1, X_2, X_3\}$. $\partial D_2$  is a simple closed curve on
$\partial{\overline{H}}$, and the Whitehead graph of $\partial D_2$ with respect to $\{X_1, X_2, X_3\}$
is shown in Figure \ref{whgraphd2}. This graph is connected and has two cut vertices ${X_2}_+$ and
${X_2}^-$. By applying the Whitehead algorithm to ${X_2}_+$ we get a new graph as shown in Figure
\ref{whgraphd21}. The graph is connected and has no cut vertex. So $\partial{\overline{H}}-\partial D_2$
is incompressible in $\overline{H}$. Hence  the manifold $M_1=\overline{H}\cup (D_2\times I)$ has
incompressible boundary $S$.

\begin{figure}
\begin{center}
\includegraphics[width=4in]{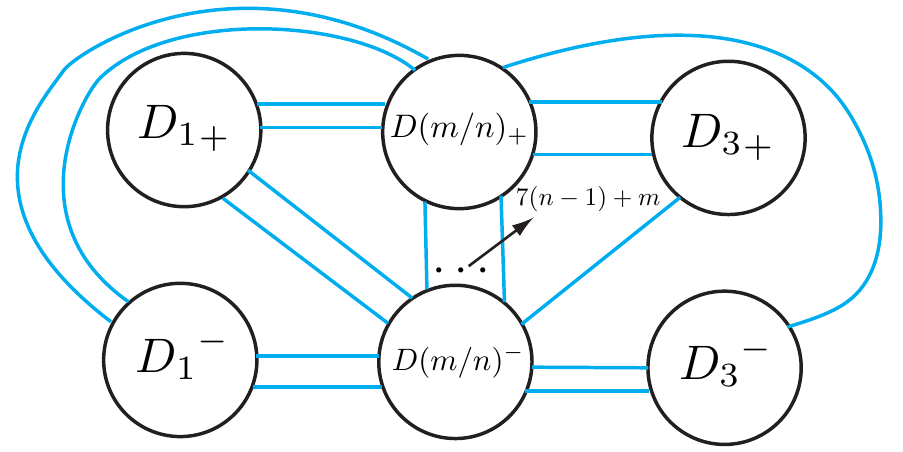}
\end{center}
\caption{\label{whgraphx4} The Whitehead graph of $\partial X$ with respect to $\{D_1, D_3, D(m/n)$.}
\end{figure}

On the other hand, compress $C(m/n)$ along $D_2$, we get a handlebody $\overline{C}$ of genus three with
$\{ D_1, D_3, D(m/n)\}$ as a  disk system. $\partial X$ is a simple closed curve on
$\partial{\overline{C}}$. The Whitehead graph of $\partial X$ with respect of $\{ D_1, D_3, D(m/n)\}$ is
shown in Figure \ref{whgraphx4}. The graph is connected and has no cut vertex. So
$\partial{\overline{C}}-\partial X$ is incompressible in $\overline{C}$. So the manifold
$M_2=\overline{C}\cup (X\times I)$ has incompressible boundary.

Thus  $S=\partial M_1=\partial M_2$  is incompressible in $M_K(m/n)=H'\cup C(m/n)$. As $S$ is contained
in $M_K$, it is also an essential surface in $M_K$.

In general, for a pretzel knot $K=(p_1, p_2, \cdots, p_k)$ as in Case 2, the proof is similar.

A regular neighborhood, $H$, is a handle body of genus $k$. We deform $H$, such that it's exterior $H'$
is a standard handlebody in $\mathbb{S}^3$. The knot complement $M_K=\mathbb{S}^3\setminus K$ has a
Heegaard splitting $M_K=H'\cup C$, where $C$ is a compression body obtained by attaching $k-1$ $1$-handles
to the positive boundary $\partial{M_K}\times[1]$. $H'$ has a meridian disk system $\{X, X_1, \cdots,
X_{k-1}\}$, and $C$ has a meridian disk system $\{D_1, \cdots, D_{k-1}\}$.

\begin{figure}
\begin{center}
\includegraphics[width=5.8in]{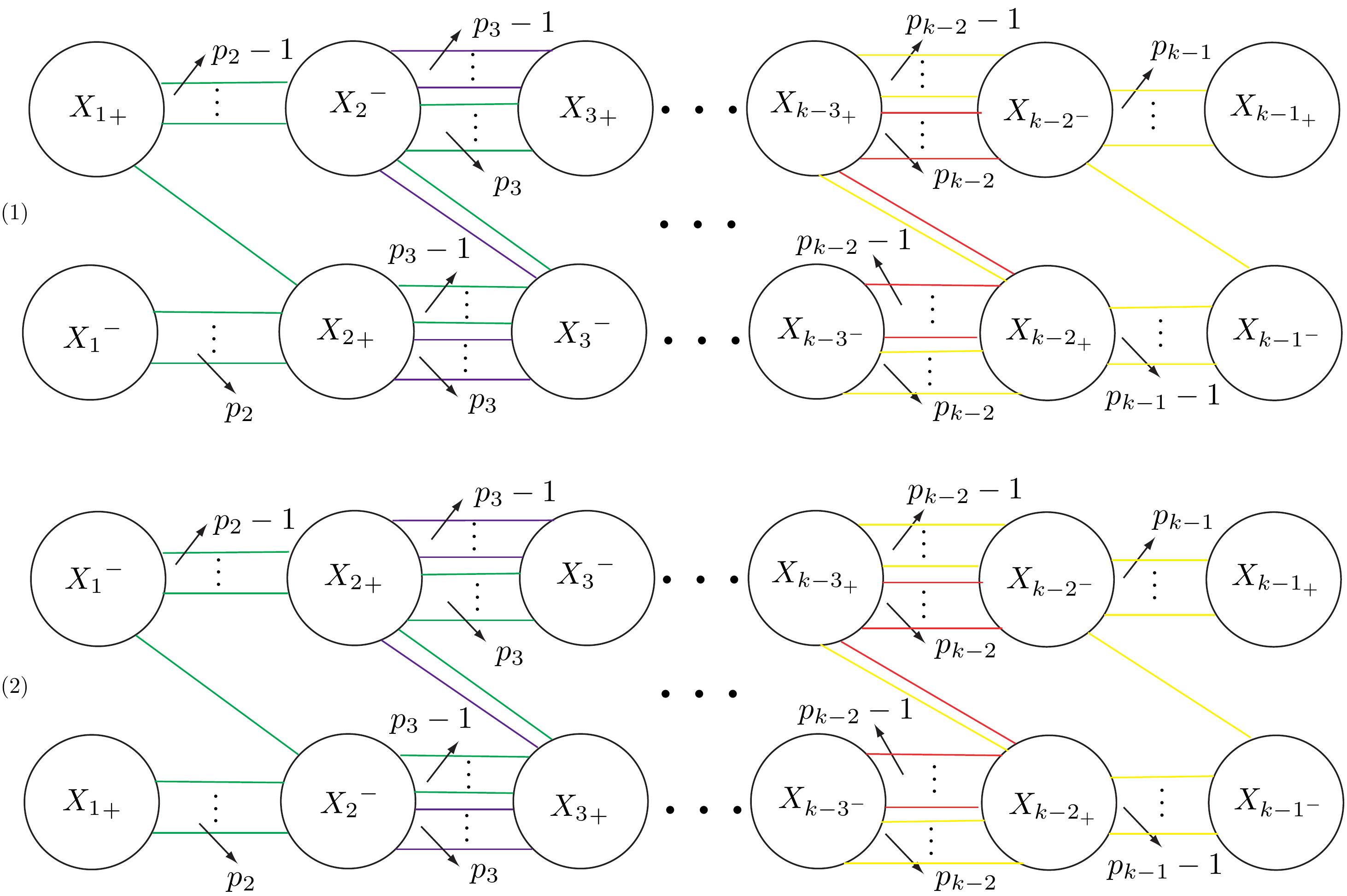}
\end{center}
\caption{\label{whgraphd2dk4} The Whitehead graph of $\{\partial D_2, \cdots, \partial D_{k-2}\}$ with
respect to $\{X_1, \cdots, X_{k-1}\}$.}
\end{figure}

Since $\{X\}$ is disjoint from $\{D_2, \cdots, D_{k-2}\}$, this Heegaard splitting is weakly reducible.
By a similar argument, we can show that the genus two surface $S$ obtained by compressing $\partial H'$
using $X$ and $D_2, \cdots, D_{k-2}$ is essential in $M_K$. Moreover, we can also show that $S$ remains
incompressible in the manifold $M_K(m/n)$ for all $m/n\ne 1/0$. In fact if $\overline{H}$ is the genus
$k-1$ handlebody obtained by compressing $H'$ along $X$, it has a   disk system $\{X_1, \cdots,
X_{k-1}\}$. $\{\partial D_2, \cdots, \partial D_{k-2}\}$ is a family of simple closed curves on
$\partial{\overline{H}}$, and the Whitehead graph of $\{\partial D_2, \cdots, \partial D_{k-2}\}$ with
respect to $\{X_1, \cdots, X_{k-1}\}$ is shown in Figure \ref{whgraphd2dk4}. 
If $k$ is even, the graph
looks like (1), if $k$ is odd, the graph looks like (2). We can check that this graph satisfies  all the
conditions of the multi-handle addition theorem. hence the manifold $M_1=\overline{H}\cup (D_2\times
I)\cup\cdots\cup(D_{k-2}\times I)$ has incompressible boundary $S$.

\begin{figure}
\begin{center}
\includegraphics[width=5in]{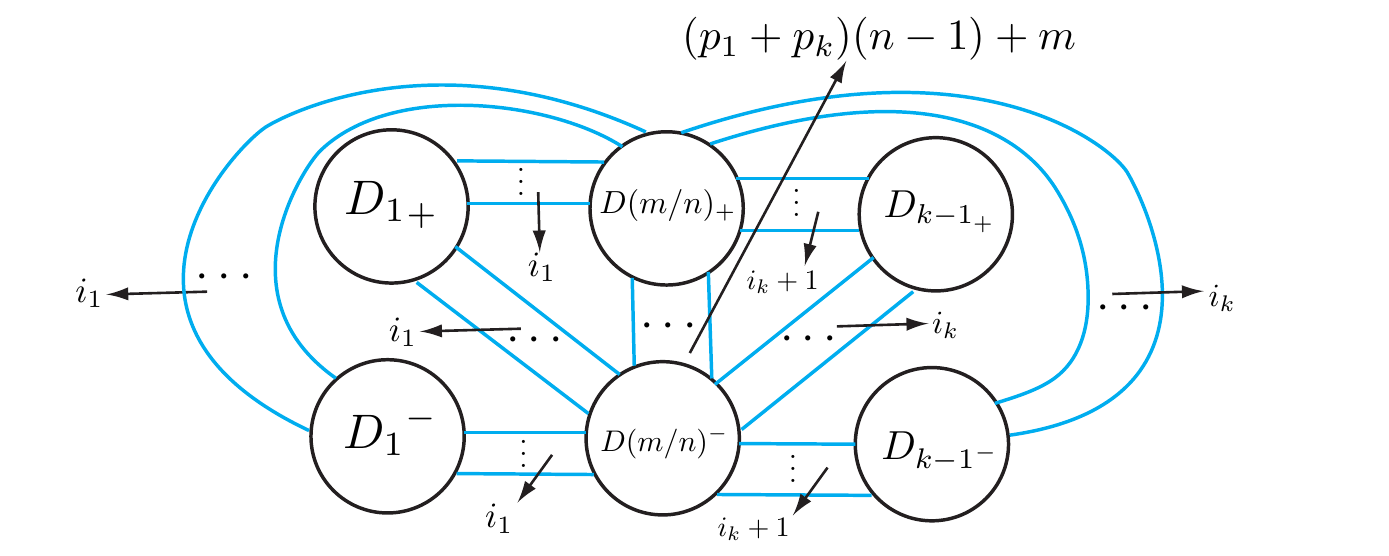}
\end{center}
\caption{\label{whgraphxg4} The Whitehead graph of $\partial X$ with respect to $\{D_1, D_{k-1},
D(m/n)$.}
\end{figure}

On the other hand, let $C(m/n)$ be the handlebody obtained by Dehn filling $C$ with slope $m/n$ and let
$D(m/n)$ be a meridian disk of the filling torus. Compress $C(m/n)$ along $D_2, \cdots, D_{k-2}$, we get
a handlebody $\overline{C}$ of genus three, and $\{ D_1, D_{k-1}, D(m/n)\}$ gives  a  disk system.
$\partial X$ is a simple closed curve on $\partial{\overline{C}}$ whose  Whitehead graph with respect of
$\{ D_1, D_{k-1}, D(m/n)\}$ is shown in Figure \ref{whgraphxg4}. The graph is connected and has no cut
vertex. So $\partial{\overline{C}}-\partial X$ is incompressible in $\overline{C}$. Thus  the manifold
$M_2=\overline{C}\cup (X\times I)$ has incompressible boundary $S$.

Now we have $S=\partial M_1=\partial M_2$, and so $S$ is incompressible in $M_K(m/n)=H'\cup C(m/n)$.
Notice that $S$ is contained in $M_K$, so $S$ is also an essential surface in $M_K$.

We finished the proof of Proposition \ref{prop}.

\end{document}